\magnification=\magstep1
\hsize=16.5 true cm 
\vsize=23.6 true cm
\font\bff=cmbx10 scaled \magstep1
\font\bfff=cmbx10 scaled \magstep2
\font\bffg=cmbx10 scaled \magstep3
\font\bffgg=cmbx10 scaled \magstep5
\parindent0cm
\font\boldmas=msbm10           %
\def\Bbb#1{\hbox{\boldmas #1}} %

\centerline{\bffgg Counting overweight spaces}
\bigskip
\centerline{\bfff Gerald Kuba}
\bigskip\medskip\smallskip
{\bff 1. Introduction}
\medskip
Write $\,|M|\,$ for the cardinal number (the {\it size}) of a set $\,M\,$
and define $\,{\bf c}:=|{\Bbb R}|=2^{\aleph_0}\,$.
We use $\,\kappa,\lambda,\mu\,$ throughout to stand for 
{\it infinite} cardinal numbers.
As usual, $\,w(X)\,$ denotes the weight of a
topological space $\,X\,$.
Naturally, $\,w(X)\leq 2^{|X|}\;$ and $\;|X|\leq 2^{w(X)}\;$ 
for every infinite T$_0$-space $\,X\,$. 
It is trivial that $\,w(X)\leq |X|\,$ for every infinite, 
first countable space $\,X\,$  and well-known (see [2, 3.3.6]) that
$\,w(X)\leq |X|\,$ 
for every compact \hbox{Hausdorff space $\,X\,$.}
Furthermore, $\,w(X)\geq |X|\,$
for every infinite, scattered T$_0$-space $\,X\,$ (see Lemma 1 below).
\medskip
According to the title, we are concerned with 
topological spaces $\,X\,$ satisfying the strict 
inequality $\,w(X)>|X|\,$. While the extreme case $\,w(X)=2^{|X|}\,$ 
is of natural interest, to investigate the case
$\;|X|<w(X)<2^{|X|}\;$ is reasonable in view of the following
remarkable fact.
\medskip
(1.1)$\;$ {\it It is consistent with {\rm ZFC} set theory that $\,\mu<\lambda\,$
implies $\,2^\mu<2^\lambda\,$ and that for every 
regular $\,\kappa\,$ there exist precisely
$\,2^\kappa\,$ cardinals $\,\lambda\,$ with
$\;\kappa<\lambda<2^\kappa\,$.}
\medskip
(A short explanation why (1.1) is true is given in Section 2.)
For fundamental enumeration theorems 
about spaces $\,X\,$ with $\,w(X)\leq|X|\,$
see [3], [5], [6], [7], [8], [9], [11].
However, it would be artificial to avoid an overlap 
with these enumeration theorems and hence in the following we include the 
case $\,w(X)=|X|\,$. The benefit of this inclusion is that we will also
establish several new enumeration theorems about spaces $\,X\,$
with $\,w(X)=|X|\,$.
A short proof of the following basic estimate
is given in the next section.
\medskip
(1.2)$\;$ {\it If $\,\theta\,$ is an infinite cardinal and
$\,{\cal F}\,$ is a family of mutually non-homeomorphic
infinite {\rm T}$_0$-spaces such that 
$\;\max\{|X|,w(X)\}\leq \theta\;$
for every $\,X\in{\cal F}\,$ then $\;|{\cal F}|\leq 2^\theta\,$.}
\smallskip\medskip
For abbreviation let us call a Hausdorff space $\,X\,$ {\it almost discrete}
if and only if $\,X\setminus\{x\}\,$
is a discrete subspace of $\,X\,$ for some $\,x\in X\,$. 
Recall that a space is {\it perfectly normal} when it is normal 
and every closed set is a G$_\delta$-set.
Note that every subspace of a perfectly normal space is 
perfectly normal. 
Recall that a normal space 
is {\it strongly zero-dimensional} if and only if for every closed set $\,A\,$
and every open set $\,U\supset A\,$ there is 
an open-closed set $\,V\,$ with $\;A\subset V\subset U\,$.
Our first goal is to prove the following enumeration theorem. 
\medskip\smallskip 
{\bf Theorem 1.} {\it If $\;\kappa\leq \lambda\leq 2^\kappa\;$ then there exist 
$\,2^{\lambda}\,$ mutually non-homeomorphic scattered,  strongly 
zero-dimensional, hereditarily paracompact, perfectly 
normal spaces $\,X\,$ with 
$\;|X|=\kappa\;$ and $\;w(X)=\lambda\,$. 
In case that $\,\lambda\leq 2^\mu<2^\lambda\,$ for some $\,\mu\,$
it can be accomplished that all these spaces are also almost 
discrete. Moreover, it can be accomplished that all these spaces
are almost discrete and extremally disconnected in case that 
$\,\lambda=2^\mu\,$ for some $\,\mu\,$ (which includes the case 
$\,\lambda=2^\kappa\,$).} \medskip 
\medskip          
Since every scattered Hausdorff space is totally disconnected,
the following theorem is a noteworthy
counter\-part of Theorem 1.
For abbreviation, let us call a space $\,X\,$ {\it almost metrizable}
if and only if $\,X\,$ is perfectly normal and $\,X\setminus\{x\}\,$
is metrizable for some $\,x\in X\,$. 
In view of Lemma 3 in Section 3, almost metrizable space
are hereditarily paracompact. 
\medskip
{\bf Theorem 2.} {\it If 
$\;{\bf c}\leq \kappa\leq \lambda\leq 2^\kappa\;$ 
then there exist $\,2^{\lambda}\,$ mutually non-homeomorphic 
pathwise connected, locally pathwise connected, almost metrizable
spaces of size 
$\,\kappa\,$ and weight $\,\lambda\,$.}
\medskip
The restriction $\,{\bf c}\leq \kappa\,$ in Theorem 2 is inevitable 
because if $\,X\,$ is an infinite, {\it pathwise} connected Hausdorff space 
then $\,X\,$ is {\it arcwise} connected (see [2, 6.3.12.a])
and hence $\,{\bf c}=|[0,1]|\leq |X|\,$.
However, for infinite, {\it connected} Hausdorff spaces $\,X\,$ 
the restriction $\,{\bf c}\leq |X|\,$ is not justified 
and we can prove the following theorem. Note that, by applying (1.1) for
$\,\kappa=\aleph_0\,$, the existence of $\,{\bf c}\,$ infinite cardinals 
$\,\kappa<{\bf c}\,$ is consistent with ZFC. 
\medskip
{\bf Theorem 3.} {\it If $\,\kappa<{\bf c}\,$ and  
$\;\kappa\leq \lambda\leq 2^\kappa\;$ 
then there exist $\,2^{\lambda}\,$ mutually non-homeomorphic 
connected and locally connected Hausdorff spaces 
of size $\,\kappa\,$ and weight $\,\lambda\,$.
In particular, up to homeomorphism
there exist precisely $\,2^{\bf c}\,$ 
countably infinite, connected, locally connected Hausdorff spaces 
and precisely $\,{\bf c}\,$ 
countably infinite, connected, locally connected, second countable 
Hausdorff spaces.}
\medskip
No space provided by Theorem 3 is completely regular
because, naturally, 
every completely regular space of size smaller than $\,{\bf c}\,$
and greater than $\,1\,$
is totally disconnected. 
Moreover, every countably infinite, regular space is totally disconnected
(see [2, 6.2.8]). The {\it connected} spaces provided by Theorem 3 are 
{\it totally pathwise disconnected} since they are Hausdorff spaces
of size smaller than $\,{\bf c}\,$. Therefore the following 
counterpart of Theorem 2 is worth mentioning. 
\medskip
{\bf Theorem 4.} {\it If 
$\;{\bf c}\leq \kappa\leq \lambda\leq 2^\kappa\;$ 
then there exist $\,2^{\lambda}\,$ mutually non-homeomorphic 
connected, totally pathwise disconnected, nowhere locally connected,
almost metrizable spaces of size $\,\kappa\,$ and weight 
$\,\lambda\,$.}                   
\bigskip
{\bff 2. Some explanations and preparations}
\medskip
Referring to Jech's profound text book [4],
a proof of (1.1) can be carried out as follows.
Define in G\"odel's universe L for every regular cardinal $\,\kappa\,$
a cardinal number $\,\theta(\kappa)\,$ by 
$\;\theta(\kappa)\,:=\,\min\{\,\mu\;|\;\mu=\aleph_\mu \;\land\;
{\rm cf}\,\mu=\kappa^+\,\}\,$.
Then $\;|\{\,\lambda\;|\;\kappa<\lambda<\theta(\kappa)\,\}|=\theta(\kappa)\;$ 
holds in every generic extension of L.
By applying Easton's theorem [4, 15.18] one can 
create an Easton universe E generically extending L
such that the continuum function $\;\kappa\mapsto 2^\kappa=\kappa^+\;$
in L is changed into $\;\kappa\mapsto 2^{\kappa}=g(\kappa)\;$ in E
with $\,g(\kappa)=\theta(\kappa)\,$ for every regular cardinal $\,\kappa\,$.
So in E we have $\;|\{\,\lambda\;|\;\kappa<\lambda<2^\kappa\,\}|=2^\kappa\;$ 
for every regular $\,\kappa\,$. By definition, in E we have 
$\,2^\alpha<2^\beta\,$ whenever $\,\alpha,\beta\,$ are regular cardinals
with $\,\alpha<\beta\,$. Therefore and in view of 
[4] Theorem 5.22 and [4] Exercise 15.12,
if $\,\mu\,$ is singular in E
then $\,2^\mu\,$ is a successor cardinal in E while 
$\,2^\kappa\,$ is a limit cardinal in E for every regular $\,\kappa\,$
in E. Consequently, in E we have 
$\,2^\mu<2^\lambda\,$ whenever $\,\mu,\lambda\,$ are arbitrary cardinals
with $\,\mu<\lambda\,$. 
\medskip
In order to verify (1.2), first of all it is clear 
that a topological space 
$\,(X,\tau)\,$ has a basis of size $\,\lambda\leq|\tau|\,$
if and only if $\;w(X)\leq\lambda\,$.
Let $\,S\,$ be an infinite set of size $\,\nu\,$ and let $\,P\,$
be the power set of $\,S\,$, whence $\,|P|=2^\nu\,$.
Let $\,\mu(\nu,\lambda)\,$ denote the 
total number of all topologies $\,\tau\,$
on $\,S\,$ such that $\,(S,\tau)\,$
has a basis $\,B\,$ of size $\,\lambda\,$. 
Clearly, $\,\mu(\nu,\lambda)=0\,$ if $\,\lambda>2^\nu\,$. 
For $\,\lambda\leq 2^\nu\,$ we have 
$\;\mu(\nu,\lambda)\leq |P|^\lambda=
\max\{2^{\nu},2^{\lambda}\}\,$.
So if $\,\theta\,$ and $\,{\cal F}\,$ satisfy the assumption in (1.2)
then $\,|{\cal F}|\,$ is not greater than the sum $\,\Sigma\,$ of all
cardinals $\;\mu(\nu,\lambda)\;$ with $\,(\nu,\lambda)\,$ running through 
the set $\;Q\,:=\,\{\,\kappa\;|\;\kappa\leq\theta\,\}^2\,$.
Thus from $\,\mu(\nu,\lambda)\leq 2^\theta\,$ for all 
$\,(\nu,\lambda)\in Q\,$ we derive
$\,\Sigma\leq 2^\theta\,$ and this concludes the proof of (1.2).

\medskip
In the following we write down a short proof 
of an important fact mentioned in the previous section.
\smallskip
{\bf Lemma 1.} {\it If $\,X\,$ is an infinite scattered {\rm T}$_0$-space
then $\,w(X)\geq|X|\,$.}
\smallskip
{\it Proof.} Since $\,X\,$ is infinite and T$_0$, no basis of $\,X\,$ is 
finite. Assume that $\,\lambda:=w(X)<|X|\,$ and
let $\,{\cal B}\,$ be a basis of $\,X\,$ 
with $\;|{\cal B}|=\lambda\,$. Let $\,X^*\,$ denote the set 
of all $\,x\in X\,$ such that $\;|U|>\lambda\;$ for every 
neighborhood $\,U\,$ of $\,x\,$. Then
$\;\;X\setminus X^*\;\subset\;
\bigcup\,\{\,U\in{\cal B}\;\,|\,\;|U|\leq \lambda\,\}\;\;$  
and hence $\;|X\setminus X^*|\leq \lambda\,$. 
Consequently, $\,X^*\not=\emptyset\,$ 
and if $\,x\in X^*\,$ and $\,U\,$  
is a neighborhood of $\,x\,$ 
then $\;|X^*\cap U|>\lambda\;$ (since $\,|U|>\lambda\,$).
Therefore, the nonempty set $\,X^*\,$ is dense in itself
and hence the space $\,X\,$ is not scattered, {\it q.e.d.}

\medskip
In order to settle the case $\,2^\kappa=2^\lambda\,$
in Theorems 1 and 2 and 4 we 
will apply the following two enumeration theorems 
about metrizable spaces. 
Note that, other than in the model E which proves (1.1), 
for $\;\kappa<\lambda\leq 2^\kappa\;$
we can rule out $\,2^\kappa=2^\lambda\,$
only in case that $\,\lambda=2^\kappa\,$.
(Thus the following two propositions can be ignored
if Theorems 1 and 2 and 4 are only read as enumeration theorems 
about spaces $\,X\,$ of 
{\it maximal possible weights} $\,2^{|X|}\,$.)
\medskip
Let $\,X+Y\,$ denote the topological sum of two Hausdorff spaces
$\,X\,$ and $\,Y\,$. (So $\,X+Y\,$ is a space $\,S\,$ 
such that $\,S=\tilde X\cup\tilde Y\,$ for disjoint
open subspaces $\,\tilde X,\tilde Y\,$ of $\,S\,$
where $\,\tilde X\,$ is homeomorphic to $\,X\,$ and 
$\,\tilde Y\,$ is homeomorphic to $\,Y\,$.)
If $\,Y=\emptyset\,$ then we put $\,X+Y=X\,$.
\medskip
{\bf Proposition 1.} {\it For every $\,\kappa\,$ 
there is a family $\;{\cal H}_\kappa\;$ 
of mutually non-homeomorphic scattered, strongly zero-dimensional
metrizable spaces
of size $\,\kappa\,$ such that 
$\,|{\cal H}_\kappa|=2^\kappa\,$ and if $\,D\,$ is any discrete 
space (including the case $\,D=\emptyset\,$) then 
the spaces $\;H_1+D\;$ and $\;H_2+D\;$ are never 
homeomorphic for distinct $\,H_1,H_2\in{\cal H}_\kappa\,$.}
\medskip
By Lemma 1 and since $\,w(Y)\leq|Y|\,$
for every metrizable space $\,Y\,$, 
we have $\,w(X)=|X|\,$ 
for every $\,X\in{\cal H}_\kappa\,$.
Proposition 1 can be verified by considering the 
spaces constructed in [9] which prove [9] Theorem 1. 
Because these spaces $\,X\,$ are revealed as  mutually
non-homeomorphic ones by investigating the $\alpha$th Cantor derivative
$\,X^{(\alpha)}\,$ for every ordinal $\,\alpha>0\,$. And, naturally, if $\,X\,$
is any space and $\,D\,$ is discrete then  
$\;(X+D)^{(\alpha)}=X^{(\alpha)}\;$ for every $\,\alpha>0\,$.
The following proposition is proved in [6] Section 4.
\medskip
{\bf Proposition 2.} {\it For every $\,\kappa\geq {\bf c}\,$
there is a family $\;{\cal P}_\kappa\;$ 
of mutually non-homeomorphic pathwise connected, locally pathwise
connected, complete metric spaces of size and weight $\,\kappa\,$
such that $\,|{\cal P}_\kappa|=2^\kappa\,$
and if $\,H\in {\cal P}_\kappa\,$ then $\,H\,$ contains a noncut point
and the cut points of $\,H\,$ lie dense in $\,H\,$.}
\bigskip
{\bff 3. Almost discrete and almost metrizable spaces}
\medskip
In accordance with [13], a space is {\it completely normal} when 
every subspace is normal. (In [2] such spaces are called
{\it hereditarily normal}.)
\medskip
{\bf Lemma 2.} {\it If $\,X\,$ is a Hausdorff space and $\,z\in X\,$ 
such that $\;X\setminus \{z\}\;$ is a discrete subspace of $\,X\,$ 
then $\,X\,$ is scattered and completely normal 
and strongly zero-dimensional.}
\smallskip
{\it Proof.} Put $\;Y:=X\setminus\{z\}\,$.
Since $\,Y\,$ is a discrete and open subspace of $\,X\,$,
every nonempty subset of $\,X\,$ contains an isolated point, whence
$\,X\,$ is scattered. 
Let $\;A,B\subset X\;$ with
$\;\overline A \cap B\,=\,A\cap\overline B\,=\,\emptyset\,$.
If $\;z\,\not\in\,A\cup B\;$ then $\;A,B\,\subset \,Y\;$ and hence 
$\;A\subset U\;$ and $\;B\subset V\;$ with the two disjoint open sets
$\;U=A\;$ and $\;V=B\,$. 
Assume $\;z\,\in\,A\cup B\;$ and, say,
$\;z\in A\,$. Then $\;z\not\in \overline B\;$ and hence $\,B\subset Y\,$.
Thus $\;A\subset \tilde U\;$ and $\;B\subset \tilde V\;$ 
with the two disjoint open sets
$\;\tilde U\,=\,X\setminus \overline B\;$ and $\;\tilde V=B\,$.
So $\,X\,$ is completely normal. Finally, 
let $\,A\subset X\,$ be closed.
If $\,z\not\in A\,$ then $\,A\,$ is open. 
If $\,z\in A\,$ and $\,U\,$ is an open neighborhood of $\,A\,$
then $\,U\,$ is closed since $\;X\setminus U\,\subset\,Y\,$.
So every closed subset of $\,X\,$ 
has a basis of open-closed neighborhoods
and hence $\,X\,$ is strongly zero-dimensional, {\it q.e.d.} 
\medskip
{\bf Lemma 3.} {\it If $\,Z\,$ is a regular space such that 
$\,Z\setminus\{z\}\,$ is paracompact for some $\,z\in Z\,$ 
then $\,Z\,$ is paracompact.}
\smallskip
{\it Proof.} Let $\,{\cal U}\,$ be an open cover
of $\,Z\,$. Trivially, 
$\;{\cal U}^*\,:=\,\{\,U\setminus\{z\}\;|\;U\in{\cal U}\,\}\;$ is an open cover
of the paracompact open subspace $\;P\,=\,Z\setminus \{z\}\;$ of $\,Z\,$.
Hence we can find an open cover $\,{\cal V}^*\,$ of $\,P\,$
which is a locally finite refinement of $\,{\cal U}^*\,$.
Fix one set $\,U_z\in{\cal U}\,$ with $\,z\in U_z\,$
and choose a closed neighborhood $\,C\,$ of $\,z\,$ in the regular space 
$\,Z\,$ such that $\,C\subset U_z\,$. 
Now put 
$\;{\cal V}\;:=\;\{\,V^*\setminus C\;\,|\;\,V^*\in{\cal V}^*\,\}
\,\cup\,\{U_z\}\,.$
Clearly, $\;{\cal V}\;$ is an open cover of $\,Z\,$
which is a refinement of $\,{\cal U}\,$.
If $\,z\not=x\in Z\,$ then some neighborhood of $\,x\,$
meets only finitely many members of $\,{\cal V}^*\,$
and hence only finitely many members of $\,{\cal V}\,$.
And $\,C\,$ is 
a neighborhood of $\,z\,$ which meets $\,V\in{\cal V}\,$
if and only if $\,V=U_z\,$.
Therefore, the cover $\,{\cal V}\,$ is locally finite in $\,Z\,$
and hence $\,Z\,$ is paracompact, {\it q.e.d.}
\medskip
Since metrizability implies paracompactness
and since the union of two G$_\delta$-sets is a G$_\delta$-set,  
from Lemma 2 and Lemma 3 we derive the following two corollaries.
\medskip
{\bf Corollary 1.} {\it Let $\,X\,$ be a Hausdorff space and $\,z\in X\,$ 
such that $\;X\setminus \{z\}\;$ is a discrete subspace of $\,X\,$
and $\,\{z\}\,$ is a {\rm G}$_\delta$-set in $\,X\,$. Then 
the almost discrete space $\,X\,$ is
hereditarily paracompact 
and perfectly normal.}
\medskip
{\bf Corollary 2.} {\it Let $\,X\,$ be a regular space and $\,z\in X\,$
such that the subspace $\,X\setminus\{z\}\,$
is metrizable and $\,\{z\}\,$ is a {\rm G}$_\delta$-set in $\,X\,$.
Then $\,X\,$ is hereditarily paracompact and perfectly normal 
and hence almost metrizable.}
\bigskip
{\bff 4. The single filter topology}
\medskip
Let $\,X,z\,$ be as in Lemma 2
and consider the family $\,{\cal U}\,$ of all open neighborhoods 
of the point $\,z\,$. Since $\,\{x\}\,$ is open in $\,X\,$ 
whenever $\,z\not=x\in X\,$, the family $\,{\cal U}\,$ coincides with the 
neighborhood filter at $\,z\,$ in the space $\,X\,$.
Consequently, 
$\;{\cal U}^*\,:=\,\{\,U\setminus\{z\}\;|\;U\in{\cal U}\,\}\;$
is the power set of $\,X\setminus\{z\}\,$ if $\,z\,$ is isolated in $\,X\,$
or, equivalently, if $\,X\,$ is discrete.
And $\,{\cal U}^*\,$ is a filter on the set $\,X\setminus\{z\}\,$
if $\,z\,$ is a limit point of $\,X\,$
or, equivalently, if the discrete subspace
$\,X\setminus \{z\}\,$ is dense in $\,X\,$.
Since $\,X\,$ is Hausdorff, it is plain that $\;\bigcap{\cal U}^*=\emptyset\,$.
\medskip
Conversely, let 
$\,Y\,$ be an infinite set and $\,z\not\in Y\,$
and let $\,{\cal F}\,$ be a filter on the set $\,Y\,$.
Define a topology $\,\tau[{\cal F}]\,$ on the set $\;X:=Y\cup\{z\}\;$
by declaring $\;U\subset X\;$ open if and only if 
either $\;z\not\in U\;$ 
or $\;U=\{z\}\cup F\;$ for some $\,F\in{\cal F}\,$.
It is plain that this is a correct definition of a topology 
on the set $\,X\,$. Furthermore, $\,Y\,$ is a discrete and open and dense
subspace of $\,(X,\tau[{\cal F}])\,$, whence $\,\{z\}\,$ is closed in $\,X\,$.
It is plain that $\,(X,\tau[{\cal F}])\,$ is a Hausdorff space 
if and only if the filter $\,{\cal F}\,$ is free, 
i.e.$\,\bigcap{\cal F}=\emptyset\,$.
So by Lemma~2 the almost discrete space
$\,(X,\tau[{\cal F}])\,$ is 
hereditarily paracompact and scattered and strongly zero-dimensional
for every free filter $\,{\cal F}\,$ on $\,Y\,$.
\smallskip
For abbreviation throughout the paper let us call a filter $\,{\cal F}\,$ 
{\it $\omega$-free} if and only if 
$\,\bigcap{\cal A}=\emptyset\,$ for some countable 
$\,{\cal A}\subset{\cal F}\,$.
In view of Corollary 1 
the following statement is evident.
\medskip
(4.1) $\;$ {\it If $\,{\cal F}\,$ is a filter on $\,Y\,$ then 
$\,(X,\tau[{\cal F}])\,$ is almost discrete and perfectly normal 
if and only if $\,{\cal F}\,$ is $\omega$-free.}
\medskip
The following observation is essential for the proof of Theorem 1.
\medskip
(4.2) $\;$ {\it If $\,{\cal F}\,$ is a free filter on $\,Y\,$ then 
the almost discrete space $\,(X,\tau[{\cal F}])\,$ is extremally disconnected
if and only if $\,{\cal F}\,$ is an ultrafilter.}
\smallskip
{\it Proof.} Firstly let $\,{\cal F}\,$ be a free ultrafilter.
Let $\,U\subset X\,$ be open.
If $\,\overline U=U\,$ then $\,\overline U\,$ is open.
So assume $\,\overline U\not=U\,$. Then
$\;\overline U\,=\,U\cup\{z\}\;$ and $\,z\not\in U\,$ since 
$\,z\,$ is the only limit point in $\,X\,$.
Thus  $\,U\subset Y\,$ and  $\,z\,$ is a limit point of $\,U\,$.
Hence every open neighborhood of $\,z\,$ meets $\,U\,$.
In other words, $\;F\cap U\not=\emptyset\;$ for every 
$\,F\in{\cal F}\,$. Consequently, $\,U\in{\cal F}\,$ 
since $\,{\cal F}\,$ is an ultrafilter. Thus 
$\;\overline U\,=\,U\cup\{z\}\;$ is open in $\,X\,$,
whence  $\,(X,\tau[{\cal F}])\,$ is extremally disconnected.
Secondly, let $\,{\cal F}\,$ be a free filter and assume that
$\,(X,\tau[{\cal F}])\,$ is extremally disconnected.
Let $\,A\subset Y\,$, whence $\,A\,$ is open in $\,X\,$. 
If $\,\overline A=A\,$ then $\,X\setminus A\,$
is open and hence $\,Y\setminus A\,$ lies in $\,{\cal F}\,$.
If $\,\overline A\not=A\,$ then $\;\overline A\,=\,\{z\}\cup A\;$
is open and hence $\,A\,$ lies in $\,{\cal F}\,$.
This reveals $\,{\cal F}\,$ as an ultrafilter, {\it q.e.d.}
\medskip
{\it Remark.} If $\,|Y|=\aleph_0\,$ and 
$\,{\cal F}\,$ is a free ultrafilter on $\,Y\,$
then $\,\tau[{\cal F}]\,$ is the well-known
{\it single ultrafilter topology} (see Example 114 in [13].)
\medskip
For a filter $\,{\cal F}\,$ on $\,Y\,$ let $\,\chi({\cal F})\,$ denote the 
least possible size of a filter base which generates $\,{\cal F}\,$.
Trivially, $\;\chi({\cal F})\leq|{\cal F}|\leq 2^{|Y|}\,$.
The notation $\,\chi(\cdot)\,$ corresponds with the obvious 
fact that $\,\chi({\cal F})\,$ is the {\it character} 
of $\,z\,$ in $\,(X,\tau[{\cal F}])\,$. (The character $\,\chi(a,A)\,$
of a point $\,a\,$ in a space $\,A\,$
is the smallest possible size of a local basis at $\,a\,$
in the space $\,A\,$.)                  
Therefore, since $\,\{y\}\,$ is open in 
$\,(X,\tau[{\cal F}])\,$ for every $\,y\in Y\,$, we obtain:
\smallskip
(4.3) $\;$ {\it If $\,{\cal F}\,$ is a free filter on $\,Y\,$ then the weight 
of $\,(X,\tau[{\cal F}])\,$ is $\;\max\{|Y|,\chi({\cal F})\}\,$.}
\bigskip
{\bf Proposition 3.} {\it If $\;|Y|=\kappa\leq\lambda\leq 2^\kappa\;$ 
then there exist $\,2^{\lambda}\,$ $\omega$-free filters 
$\,{\cal F}\,$ on $\,Y\,$ such that $\;\chi({\cal F})=\lambda\,$.}
\medskip
{\it Remark.} The cardinal $\,2^{\lambda}\,$ in Proposition 3 is best 
possible. Indeed, let $\,Y\,$ be an infinite set of size $\,\kappa\,$
and let $\,\lambda\geq\kappa\,$.
Since a filter base on 
$\,Y\,$ is a subset of the power set of $\,Y\,$, there are at most 
$\,2^\lambda\,$ filter bases $\,{\cal B}\,$ on $\,Y\,$ with 
$\,|{\cal B}|=\lambda\,$. Hence $\,Y\,$ cannot carry more than 
$\,2^\lambda\,$ filters $\,{\cal F}\,$  
with $\;\chi({\cal F})=\lambda\,$.
\medskip
{\it Proof of Proposition 3.} Assume $\;|Y|=\kappa\leq\lambda\leq 2^\kappa\;$ 
and let $\,{\cal A}\,$ be a family of subsets 
of $\,Y\,$ such that $\,|{\cal A}|=2^\kappa\,$
and 
\smallskip
(4.4)\quad {\it If $\;{\cal D},{\cal E}\not=\emptyset\;$ are 
disjoint finite subfamilies of $\,{\cal A}\,$ then 
$\;\bigcap{\cal D}\,\not\subset\,\bigcup {\cal E}\,$.}\smallskip
\smallskip
A construction of such a family $\,{\cal A}\,$ is elementary, see [4, 7.7].
However, this is not enough for our purpose.
In view of the property {\it $\omega$-free}, 
we additionally have to make sure that 
the family $\,{\cal A}\,$ also contains a countably
infinite family $\,{\cal A}_\omega\,$ 
such that $\,\bigcap{\cal A}_\omega=\emptyset\,$.
By applying Lemma 8 in Section 11 for $\,\mu=\aleph_0\,$
we can assume that such a family $\,{\cal A}_\omega\subset{\cal A}\,$ exists.
\smallskip
Now put 
\smallskip
\centerline{$\;{\bf A}_\lambda\;:=\;\{\,{\cal H}\;\,|\,\;
{\cal A}_\omega\subset {\cal H}\subset {\cal A}\;
\land\;|{\cal H}|=\lambda\,\}\,.$}
\smallskip
Clearly, $\;|{\bf A}_\lambda|=(2^\kappa)^\lambda=2^\lambda\,$.
By virtue of (4.4), if for $\,{\cal H}\in {\bf A}_\lambda\,$
we put 
\smallskip
\centerline{${\cal B}_{\cal H}\;:=\;
\{\,H_1\cap\cdots\cap H_n\;
|\;n\in{\Bbb N}\;\land\;H_1,...,H_n\in{\cal H}\,\}$}
\smallskip
then $\;\emptyset\not\in{\cal B}_{\cal H}\;$ and hence 
$\,{\cal B}_{\cal H}\,$ is a filter base on $\,Y\,$.
For every $\,{\cal H}\in{\bf A}_\lambda\,$ let 
$\,{\cal F}[{\cal H}]\,$ denote the filter on $\,Y\,$ 
generated by $\,{\cal B}_{\cal H}\,$.
Clearly, $\;|{\cal B}_{\cal H}|=|{\cal H}|=\lambda\;$
for every $\,{\cal H}\in{\bf A}_\lambda\,$.
\smallskip
The filter $\,{\cal F}[{\cal H}]\,$ is $\omega$-free
because $\;{\cal A}_\omega\subset {\cal F}[{\cal H}]\;$
by definition. 
Furthermore, (4.4) implies that for distinct families
$\;{\cal H}_1,{\cal H}_2\in{\bf A}_\lambda\;$ 
the filters $\,{\cal F}[{\cal H}_1]\,$ and
$\,{\cal F}[{\cal H}_2]\,$ must be distinct.
So the family $\;\{\,{\cal F}[{\cal H}]\;|\;{\cal H}\in{\bf A}_\lambda\,\}\;$
consists of $\,2^\lambda\,$ $\omega$-free filters on $\,Y\,$.
\smallskip
It remains to verify that 
$\;\chi({\cal F}[{\cal H}])=\lambda\;$
for every $\,{\cal H}\in {\bf A}_\lambda\,$.
Assume indirectly that for some $\,{\cal H}\in {\bf A}_\lambda\,$
we have $\;\chi({\cal F}[{\cal H}])\not=\lambda\;$
and hence $\;\chi({\cal F}[{\cal H}])<\lambda\,$.
(Clearly $\;\chi({\cal F}[{\cal H}])\leq\lambda\;$
since $\;|{\cal B}_{\cal H}|=|{\cal H}|=\lambda\,$.)
Choose a filter base
$\,{\cal B}\,$ on $\,Y\,$ which 
generates the filter $\,{\cal F}[{\cal H}]\,$ such that 
$\,|{\cal B}|<\lambda\,$.
Since $\;{\cal B}\subset{\cal F}[{\cal H}]\;$
and $\,{\cal F}[{\cal H}]\,$ is generated by the filter base
$\,{\cal B}_{\cal H}\,$,
we can choose for every $\,B\in{\cal B}\,$ 
a finite set $\,{\cal H}_B\subset{\cal H}\,$
such that 
$\;B\,\supset\,\bigcap {\cal H}_B\,$.
Put $\;{\cal U}\,:=\,\bigcup_{B\in{\cal B}}{\cal H}_B\,$.
Then $\,{\cal U}\subset {\cal H}\,$
and $\,|{\cal U}|\leq|{\cal B}|<\lambda\,$.
Consequently, $\;{\cal H}\setminus {\cal U}\,\not=\,\emptyset\,$.
Choose any set $\;A\,\in\,{\cal H}\setminus {\cal U}\,$.
Then $\,A\in{\cal F}[{\cal H}]\,$ and hence
we can find a set  $\,B\in{\cal B}\,$
with $\,A\supset B\,$. Then
$\;A\,\supset\,\bigcap {\cal H}_B\;$
and hence  $\,A\in{\cal H}_B\,$ by virtue of (4.4).
But then $\,A\in{\cal U}\,$
in contradiction with choosing $\,A\,$
in $\,{\cal H}\setminus {\cal U}\,$, {\it q.e.d.}
\medskip          
Proposition 3 can be improved 
in the important case $\,\lambda=2^\kappa\,$ as follows.
\medskip
{\bf Proposition 4.} {\it On an infinite set of size 
$\,\kappa\,$ there exist precisely $\,2^{2^\kappa}\,$ $\omega$-free 
ultrafilters 
$\,{\cal F}\,$ such that $\;\chi({\cal F})=2^\kappa\,$.}
\medskip
{\it Proof.} Let $\,Y\,$ be a set of size $\,\kappa\,$.
As in the previous proof let $\,{\cal A}\,$ be a family of subsets 
of $\,Y\,$ such that $\,|{\cal A}|=2^\kappa\,$
and (4.4) holds. (Here we need not consider 
$\,{\cal A}_\omega\subset{\cal A}\,$.)
\smallskip
Let $\,{\bf A}\,$ denote the family of all subfamilies $\,{\cal G}\,$
of $\,{\cal A}\,$ such that 
$\;|{\cal G}|=2^\kappa\;$. Clearly, $\;|{\bf A}|=2^{2^\kappa}\,$.
Now for every $\,{\cal G}\in{\bf A}\,$ define 
\smallskip
\centerline{${\cal W}[{\cal G}]\;:=\;{\cal G}\,\cup\,
\{\,Y\setminus \bigcap {\cal H}\;|\;
{\cal H}\subset{\cal G}\;\land\;|{\cal H}|\geq\aleph_0\,\}\,\cup\,
\{\,Y\setminus A\;|\;A\,\in\,{\cal A}\setminus{\cal G}\,\}\,.$}
\smallskip     
A moment's reflection suffices to see that 
(4.4) implies that $\;W_1\cap\cdots\cap W_n\,\not=\,\emptyset\;$
whenever $\;W_1,...,W_n\in {\cal W}[{\cal G}]\,$. 
Hence for every $\,{\cal G}\in{\bf A}\,$ we can choose an ultrafilter 
$\,{\cal U}[{\cal G}]\,$ on $\,Y\,$ such that 
$\;{\cal U}[{\cal G}]\supset {\cal W}[{\cal G}]\,$ (see [1] 7.1).
\smallskip
If $\;{\cal G}_1,{\cal G}_2\in {\bf A}\;$ are distinct and, say, 
$\;G\,\in\,{\cal G}_1\setminus{\cal G}_2\;$ then 
$\;G\in {\cal W}[{\cal G}_1]\;$ and $\;Y\setminus G\,\in\, {\cal W}[{\cal G}_2]\;$
and hence $\;G\in {\cal U}[{\cal G}_1]\;$ and $\;G\not\in {\cal U}[{\cal G}_2]\;$
and hence the ultrafilters 
$\,{\cal U}[{\cal G}_1]\,$ and $\,{\cal U}[{\cal G}_2]\,$ are distinct as well.
Consequently, the family $\;\{\,{\cal U}[{\cal G}]\;|\;{\cal G}\in{\bf A}\,\}\;$
consists of $\,2^{2^\kappa}\,$ ultrafilters on $\,Y\,$.
All these ultrafilters are $\omega$-free because if 
$\;{\cal G}\in{\bf A}\;$ and $\,{\cal H}\,$ is a countably infinite
subset of $\,{\cal G}\,$ then by virtue of (4.4) the family 
$\;{\cal H}^*\,:=\,\{\,H\setminus\bigcap{\cal H}\;|\;H\in{\cal H}\,\}\;$
is countably infinite and it is trivial that  
$\;\bigcap{\cal H}^*=\emptyset\;$ and from $\;{\cal H}\subset {\cal W}[{\cal G}]\;$
and $\;Y\setminus\bigcap{\cal H}\,\in\,{\cal W}[{\cal G}]\;$
we derive $\;{\cal H}^*\subset {\cal U}[{\cal G}]\,$.
(Actually, 
by a deep argument from set theory
it is superfluous to verify that $\,{\cal U}[{\cal G}]\,$
is $\omega$-free, see the remark below.) 
\smallskip
Finishing the proof, we claim that 
$\;\chi({\cal U}[{\cal G}])=2^\kappa\;$ for every $\,{\cal G}\in{\bf A}\,$.
Assume indirectly that for $\,{\cal G}\in{\bf A}\,$ the ultrafilter
$\,{\cal U}[{\cal G}]\,$ is generated by a filter base 
$\,{\cal B}\,$ with $\,|{\cal B}|<2^\kappa\,$.
Since $\,{\cal G}\subset {\cal U}[{\cal G}]\,$, for every $\,G\in{\cal G}\,$
we have $\,G\supset B\,$ for some $\,B\in{\cal B}\,$.
From $\,|{\cal B}|<|{\cal G}|\,$ we derive the existence of a set 
$\,B\in{\cal B}\,$ and an infinite subset  
$\,{\cal H}\subset{\cal G}\,$ such that $\,H\supset B\,$ for every
$\,H\in{\cal H}\,$. Consequently, $\,\bigcap{\cal H}\supset B\,$
and hence $\,\bigcap{\cal H}\,\in\,{\cal U}[{\cal G}]\,$. 
This, however, is a contradiction since 
$\;Y\setminus \bigcap{\cal H}\;$
lies in $\,{\cal U}[{\cal G}]\,$ by the definition of $\,{\cal W}[{\cal G}]\,$, {\it q.e.d.}
\medskip
{\it Remark.} Our proof of Proposition 4 is elementary and purely 
set-theoretical. There is also a topological but much less
elementary way to prove Proposition~4.
First of all, if one can prove that any set of size $\,\kappa\,$
carries $\,2^{2^\kappa}\,$ ultrafilters of character $\,2^\kappa\,$
then Proposition~4 must be true. 
Because, an ultrafilter $\,{\cal F}\,$ is free if and only if 
$\,\chi({\cal F})>1\,$ and 
if a free ultrafilter $\,{\cal F}\,$
is not $\omega$-free then it is plain that $\,{\cal F}\,$ is $\sigma$-complete.
However, the existence of a $\sigma$-complete free ultrafilter
is unprovable in ZFC! (See [4, 10.2] and [4, 10.4].)
Now, consider the set $\,Y\,$ of size $\,\kappa\,$
equipped with the discrete topology
and consider the Stone-\v Cech compactification $\,\beta Y\,$ of $\,Y\,$
and its compact remainder $\,Y^*\,=\,\beta Y\setminus Y\,$. 
So the points in $\,Y^*\,$ are the free ultrafilters on $\,Y\,$
and if for $\,p\in Y^*\,$ we consider the subspace 
$\,Y\cup\{p\}\,$ of $\,\beta Y\,$ then it is clear that
the character of the ultrafilter $\,p\,$ equals 
$\,\chi(p,Y\cup\{p\})\,$.
It is a nice exercise to verify that 
$\;\chi(p,Y\cup\{p\})=\chi(p,Y^*)\;$
for every $\,p\in Y^*\,$. By embedding an appropriate 
Stone space of a Boolean algebra into $\,Y^*\,$ 
it can be proved that $\,Y^*\,$ must contain 
$\,2^{2^\kappa}\,$ points $\,p\,$ with $\;\chi(p,Y^*)=2^\kappa\,$,
see 7.13, 7.14, 7.15 in [1].
\bigskip
{\bff 5. Proof of Theorem 1}
\medskip
Assume $\;\mu\leq\kappa\leq\lambda\leq 2^\mu\;$
and let $\,Y\,$ be a set of size $\,\mu\,$.
Let $\,{\bf F}_\lambda\,$ denote a family of $\omega$-free filters 
on $\,Y\,$ such that $\;|{\bf F}_\lambda|=2^\lambda\;$
and $\,\chi({\cal F})=\lambda\,$ for every $\,{\cal F}\in{\bf F}_\lambda\,$.
Such a family exists by Proposition~3.
We additionally assume that if $\,\lambda=2^\mu\,$ then
every member of $\,{\bf F}_\lambda\,$
is an ultrafilter. This additional assumption is justified by 
Proposition~4. 
\medskip
Now fix $\,z\not\in Y\,$ and for every $\,{\cal F}\in{\bf F}_\lambda\,$
consider the single filter topology $\,\tau[{\cal F}]\,$
on the set $\;X\,=\,Y\cup\{z\}\;$ as in Section 4. If $\,\mu<\kappa\,$ 
then let $\,D\,$ be a discrete space of size $\,\kappa\,$.
If $\,\mu=\kappa\,$ then put $\,D=\emptyset\,$.
In both cases define the space $\,(\tilde X,\tilde\tau[{\cal F}])\,$
as the topological sum of $\,D\,$ and the space $\,(X,\tau[{\cal F}])\,$.
(So if $\,\mu=\kappa\,$ then $\,\tilde X=X\,$ and 
$\,\tilde\tau[{\cal F}]=\tau[{\cal F}]\,$.)
Clearly, $\,\tilde X\,$ 
is almost discrete, scattered, 
strongly zero-dimensional, hereditarily paracompact, and
perfectly normal. Furthermore, $\,w(\tilde X)=\lambda\,$ and 
$\,|\tilde X|=\kappa\,$.
If $\,\lambda=2^\mu\,$ 
then the space $\,\tilde X\,$ is also 
extremally disconnected by virtue of (4.2).
\medskip
Obviously, $\,\tilde\tau[{\cal F}_1]\not=\tilde\tau[{\cal F}_2]\,$
whenever the filters 
$\,{\cal F}_1,{\cal F}_2\in{\bf F}_\lambda\,$ are distinct.
(For if $\,{\cal F}_1,{\cal F}_2\in{\bf F}_\lambda\,$ and 
$\;F\in{\cal F}_1\setminus{\cal F}_2\;$ then 
$\,F\cup\{z\}\,$ is $\tilde\tau[{\cal F}_1]$-open 
but not $\tilde\tau[{\cal F}_2]$-open.) 
Consequently, the family 
$\;{\cal T}_\lambda\,:=\,
\{\,\tilde\tau[{\cal F}]\;|\;{\cal F}\in{\bf F}_\lambda\,\}\;$
is of size $\,2^\lambda\,$.
\medskip
We distinguish the two cases $\,2^\lambda>2^\mu\,$ and 
$\,2^\lambda\leq 2^\mu\,$. Assume firstly that $\,2^\lambda>2^\mu\,$
or, equivalently, that $\;|{\cal T}_\lambda|>2^\mu\,$.
Define an equivalence relation $\,\sim\,$ on $\,{\cal T}_\lambda\,$
by $\,\tau_1\sim \tau_2\;$ if and only if
the spaces $\,(\tilde X,\tau_1)\,$
and $\,(\tilde X,\tau_2)\,$ are homeomorphic.
We claim that the size of an equivalence class cannot be greater
than $\,2^\mu\,$. 
\smallskip
This is clearly true if $\,\mu=\kappa\,$
because there are only $\,2^\mu\,$ permutations on $\,X\,$.
So assume $\,\mu<\kappa\,$. If $\,\tau\in {\cal T}_\lambda\,$
then in the space $\,(\tilde X,\tau)\,$ the point $\,z\,$ 
is the only limit point and every neighborhood $\,U\,$
of $\,z\,$ is open-closed.
As a consequence, for $\,\tau_1,\tau_2\in{\cal T}_\lambda\,$ 
the spaces $\,(\tilde X,\tau_1)\,$ and $\,(\tilde X,\tau_2)\,$ 
are homeomorphic if and only if there 
is a homeomorphism $\,\varphi\,$ from the $\tau_1$-subspace 
$\,X\,$ of $\,\tilde X\,$ onto some $\tau_2$-open-closed 
subspace of $\,\tilde X\,$. 
Indeed, if $\,f\,$ is a homeomorphism from 
$\,(\tilde X,\tau_1)\,$ onto $\,(\tilde X,\tau_2)\,$
then put 
$\,\varphi(x)=f(x)\,$ for every $\,x\in X\,$ and $\,\varphi\,$ 
fits since $\,f(z)=z\,$.
Conversely, if $\,\varphi\,$ 
is a homeomorphism from the $\tau_1$-subspace 
$\,X\,$ of $\,\tilde X\,$ onto some $\tau_2$-open-closed 
subspace of $\,\tilde X\,$
and $\,g\,$ is any bijection from $\;\tilde X\setminus X\;$ onto
$\;\tilde X\setminus \varphi(X)\;$ then it is plain that 
a homeomorphism $\,f\,$ from 
$\,(\tilde X,\tau_1)\,$ onto $\,(\tilde X,\tau_2)\,$
is defined by $\;f(x)=\varphi(x)\;$ for $\,x\in X\,$
and $\;f(x)=g(x)\;$ for $\,x\not\in X\,$.
(Note that $\;|\tilde X\setminus X|=|\tilde X\setminus \varphi(X)|\;$
since $\,\mu<\kappa\,$.)
Therefore, since 
there are precisely $\,\kappa^\mu\,$ mappings from $\,X\,$ into $\,\tilde X\,$,
the size of an eqivalence class in $\,{\cal T}_\lambda\,$
cannot exceed $\,\kappa^\mu\,$. And from 
$\;2<\mu\leq\kappa\leq 2^\mu\;$ we derive 
$\;2^\mu\leq \mu^\mu\leq \kappa^\mu\leq (2^\mu)^\mu=2^\mu\;$
and hence $\,\kappa^\mu=2^\mu\,$.
\smallskip
So the size of an equivalence class can indeed not be greater
than $\,2^\mu\,$. Consequently, $\;|{\cal T}_\lambda|>2^\mu\;$
implies that the total number of all equivalence classes equals 
$\,|{\cal T}_\lambda|=2^{\lambda}\,$.  
Thus by choosing one topology 
in each equivalence class we obtain $\,2^{\lambda}\,$
mutually non-equivalent topologies $\,\tau\in {\cal T}_\lambda\,$
and hence the $\,2^{\lambda}\,$ corresponding spaces $\,(\tilde X,\tau)\,$ 
are mutually non-homeomorphic. This settles 
the case $\,2^\lambda>2^\mu\,$.
In particular, we have 
already proved the second and the third statement 
in Theorem 1 because, under the assumption 
$\,\kappa\leq\lambda\leq 2^\kappa\,$,
if $\,\lambda=2^\mu\,$ for some $\,\mu\,$ then 
$\,\lambda=2^\mu\,$ (and hence $\,2^\lambda>2^\mu\,$)  
for some $\,\mu\leq\kappa\,$
and if $\,\lambda\leq 2^\mu<2^\lambda\,$ and $\,\mu>\kappa\,$
then $\,2^\kappa\leq 2^\mu<2^\lambda\,$
and hence $\,2^\lambda>2^{\mu'}\,$ for $\,\mu'=\kappa\,$.
\medskip
Secondly assume that $\;2^\lambda\leq 2^\mu\,$.
Then we have $\;2^\lambda=2^\kappa\;$ 
since $\;\mu\leq\kappa\leq\lambda\;$ implies 
$\;2^\mu\leq 2^\kappa\leq 2^\lambda\,$.
So in order to conclude the proof of Theorem 1 
we assume $\;\kappa\leq\lambda\leq 2^\kappa=2^\lambda\,$.
(Then, of course, $\;\kappa\leq\lambda< 2^\kappa=2^\lambda\,$.)
Since the special case $\,\kappa=\lambda\,$ is 
settled by Proposition 1, we also assume $\,\kappa<\lambda\,$.
For two spaces $\,X_1\,$ and $\,X_2\,$ let, again,
$\,X_1+X_2\,$ denote the topological sum of $\,X_1\,$ and $\,X_2\,$. 
Let $\,{\cal H}_\kappa\,$ be a family provided by Proposition 1.
Due to metrizability, every space in $\,{\cal H}_\kappa\,$
is perfectly normal and hereditarily paracompact.
\medskip
By considering an appropriate
single filter topology on a set of size $\,\kappa\,$,
we can choose a perfectly normal space $\,Z\,$ of size $\,\kappa\,$ 
such that for some point $\,z\in Z\,$
the subspace $\,Z\setminus\{z\}\,$ is discrete
and $\;\chi(z,Z)=\lambda\,$. (Consequently, $\,w(Z)=\lambda\,$.)
For every space $\,H\in{\cal H}_\kappa\,$ 
consider the topological sum $\;H+Z\,$.
Of course, the topological sum of two paracompact spaces is paracompact
and $\;(H+Z)\setminus\{z\}\,=\,H+(Z\setminus\{z\})\;$
for every $\,H\in{\cal H}_\kappa\,$.
Consequently, for every $\,H\in{\cal H}_\kappa\,$
the space 
$\,H+Z\,$ is scattered and strongly zero-dimensional and perfectly normal
and hereditarily paracompact and $\;|H+Z|=|H|=\kappa\;$ and  
$\;w(H+Z)=\max\{w(H),w(Z)\}=\max\{\kappa,\lambda\}=\lambda\,$.
Therefore, since $\,|{\cal H}_\kappa|=2^\kappa\,$, the case 
$\,2^\lambda=2^\kappa\,$ in Theorem 1 is settled 
by showing that for two distinct (and hence non-homeomorphic)
metrizable spaces $\;H_1,H_2\in{\cal H}_\kappa\;$ 
the two spaces $\,H_1+Z\,$ and $\,H_2+Z\,$ are never homeomorphic.
Assume that 
$\;H_1,H_2\in{\cal H}_\kappa\;$ and that $\,f\,$ is a homeomorphism 
from $\,H_1+Z\,$ onto $\,H_2+Z\,$. Then $\,f(z)=z\,$
since $\,w((H_i+Z)\setminus\{z\})=\kappa<\lambda\,$
and $\;\chi(z,H_i+Z)=\chi(z,Z)=\lambda\,$.
Consequently, $\,f\,$ maps $\,(H_1+Z)\setminus\{z\}\,$
onto $\,(H_2+Z)\setminus\{z\}\,$. Therefore, since 
$\,Z\setminus\{z\}\,$ is discrete and
$\;(H+Z)\setminus\{z\}\,=\,H+(Z\setminus\{z\})\;$
for every $\,H\in{\cal H}_\kappa\,$,
we have $\,H_1=H_2\,$ in view of Proposition 1.
\bigskip
{\bff 6. Proof of Theorem 2}
\medskip
In order to find a natural way to prove
Theorem 2 (and also Theorem 4) we give a short proof
of the following consequence of Theorem 2.
\medskip
(6.1) {\it If $\;{\bf c}\leq \kappa\leq\lambda\leq 2^\kappa\;$ 
then there exist $\,2^{\lambda}\,$ mutually non-homeomorphic 
pathwise connected, paracompact Hausdorff spaces of size 
$\,\kappa\,$ and weight $\,\lambda\,$.}
\medskip 
(6.1) can easily be derived from Theorem 1 as follows.
Assume $\,{\bf c}\leq \kappa\leq\lambda\leq 2^\kappa\,$. 
By Theorem 1 there exists a family $\,{\cal P}\,$
of $\,2^\lambda\,$ mutually non-homeomorphic, totally disconnected,
paracompact Hausdorff spaces 
$\,X\,$ of size $\,\kappa\,$ and weight $\,\lambda\,$.
For every $\,X\in{\cal F}\,$ let $\,{\cal Q}(X)\,$ 
denote the quotient space of $\,X\times[0,1]\,$ by its 
closed subspace $\,X\times\{1\}\,$. The quotient space $\,{\cal Q}(X)\,$
can be directly defined as follows.
Consider the product space $\;X\times[0,1[\;$ and 
fix $\;p\not\in X\!\times\![0,1[\;$ and
put $\;{\cal Q}(X)\,:=\,\{p\}\cup(X\times[0,1[)\,$. Declare a subset $\,U\,$ of 
$\,{\cal Q}(X)\,$ open if and only if 
$\;U\setminus\{p\}\;$ is open in the product space $\;X\times[0,1[\;$ and 
$\,p\in U\,$ implies that 
$\;(U\setminus\{p\})\cup(X\times\{1\})\;$ is open 
in the space $\,X\times[0,1]\,$.
One can picture $\,{\cal Q}(X)\,$ as a {\it cone} with apex 
$\,p\,$ and all rulings $\;\{p\}\cup(\{x\}\times[0,1[)\;(x\in X)\;$ 
homeomorphic to the unit interval $\,[0,1]\,$. 
By [2, 5.1.36] and [2, 5.1.28] both 
$\,X\times[0,1]\,$ and $\,X\times[0,1[\;$ are paracompact.
Consequently, $\,{\cal Q}(X)\,$ is a regular space
and hence $\,{\cal Q}(X)\,$ is paracompact in view of Lemma 3.
It is evident that $\,{\cal Q}(X)\,$ is pathwise connected.
Trivially, $\,|{\cal Q}(X)|=\kappa\,$. 
\smallskip 
Unfortunately we can be sure that $\,w({\cal Q}(X))=\lambda\,$
for every $\,X\in{\cal P}\,$ only if $\,\lambda=2^\kappa\,$.
(Since $\,|{\cal Q}(X)|=\kappa\,$, we have $\,w({\cal Q}(X))\leq 2^\kappa\,$.
On the other hand, $\;w({\cal Q}(X))\geq w({\cal Q}(X)\setminus\{p\})=
w(X\times[0,1[)=w(X)=\lambda\,$.) The problem with the weight is
that if $\,\mu\,$ is the character of the apex $\,p\,$ then
$\;w({\cal Q}(X))=\max\{w(X\times[0,1[),\mu\}=
\max\{\lambda,\mu\}\,$. But we cannot rule out $\,\lambda<\mu\,$
if $\,\lambda<2^\kappa\,$. Of course, if 
$\,X\in{\cal P}\,$ is compact then $\,\mu=\aleph_0\,$ 
and hence $\,w({\cal Q}(X))=\lambda\,$ 
(but also $\,\lambda\leq |X|=\kappa\,$).
Fortunately, we can make the character of the apex countable 
also by harshly reducing the filter
of the neighborhoods of $\,p\,$. Let $\,{\cal Q}^*(X)\,$ 
be defined as the cone $\,{\cal Q}(X)\,$ but with the (only) difference 
that $\;U\,\subset\,\{p\}\cup(X\times[0,1[)\;$ is an open 
neighborhood of $\,p\,$ if and only if $\,U\setminus\{p\}\,$
is open in $\;X\times[0,1[\;$ and 
$\;U\,\supset\,X\times [t,1[\;$ for some $\;t\in[0,1[\,$.
Now we have $\;\chi(p,{\cal Q}^*(X))=\aleph_0\;$
and hence $\;w({\cal Q}^*(X))=w(X)\;$ for every $\,X\in{\cal P}\,$.
Of course, $\,{\cal Q}^*(X)\,$ is pathwise connected. 
By the same arguments as for $\,{\cal Q}(X)\,$, 
the space $\,{\cal Q}^*(X)\,$ is regular and paracompact. 
Finally, the 
spaces $\;{\cal Q}(X)\;(X\in{\cal P})\;$ are mutually non-homeomorphic
because every $\;X\in{\cal P}\;$ can be recovered (up to homeomorphism)
from $\,{\cal Q}(X)\,$. Indeed, since $\,X\,$ is totally disconected,
if $\,Z\,$ is the set of all $\,z\in {\cal Q}(X)\,$ such that 
$\;{\cal Q}(X)\setminus\{z\}\;$ 
remains pathwise connected then it is evident that
$\;Z\,=\,X\times\{0\}\;$ and hence $\,Z\,$ is homeomorphic with $\,X\,$.
This concludes the proof of (6.1).
\medskip
In the following proof of Theorem 2 we will also work with cones
but we cannot use the cones $\,{\cal Q}(X)\,$ or $\,{\cal Q}^*(X)\,$ 
because it is evident that if $\,X\,$ is not discrete
then neither $\,{\cal Q}(X)\,$ nor $\,{\cal Q}^*(X)\,$ is 
locally connected. Furthermore, by virtue of Corollary 2
and since $\,\{p\}\,$ is a G$_\delta$-set in the space $\,{\cal Q}^*(X)\,$,
the cone $\,{\cal Q}^*(X)\,$ is almost metrizable if and only if 
$\,X\,$ is metrizable. (But then $\,w({\cal Q}^*(X))=w(X)=\kappa\,$.)
Consequently, $\,{\cal Q}^*(X)\,$ is locally connected and 
almost metrizable if and only if $\,X\,$ is discrete. 
Now the clue in the following proof of Theorem 2 is to consider
$\,{\cal Q}^*(S)\,$ for one discrete spaces $\,S\,$ of size (and weight)
$\,\kappa\,$ and to reduce the topology of $\,{\cal Q}^*(S)\,$ 
in $\,2^\lambda\,$ ways such that the weight $\,\kappa\,$ of 
$\,{\cal Q}^*(S)\,$ is increased to $\,\lambda\,$ 
and that $\,2^\lambda\,$ non-homeomorphic spaces as desired are obtained.
First of all we need a lemma.
\medskip                                                                     
{\bf Lemma 4.} {\it If $\,n\in{\Bbb N}\,$ and $\,A\,$ is a 
topological space and $\,a\in A\,$
and $\;A_1,...,A_n\;$ are metrizable, closed subspaces of $\,A\,$
and $\;A=A_1\cup\cdots\cup A_n\;$
and $\;A_i\cap A_j=\{a\}\;$ whenever $\,1\leq i<j\leq n\,$
then the space $\,A\,$ is metrizable.}
\smallskip
{\it Proof.} Assume $\,n\geq 2\,$. Clearly, if $\,1\leq i\leq n\,$ then
$\;A_i\setminus\{a\}\,=\,
A\setminus\bigcup_{j\not=i} A_j\;$ is an open subset of $\,A\,$.
Furthermore, if $\,a\in U_i\subset A_i\,$ 
and $\,U_i\,$ is open in the subspace $\,A_i\,$ for $\,1\leq i\leq n\,$
then $\;U_1\cup\cdots\cup U_n\;$ is an open subset of the space $\,A\,$.
(Because if $\,V_i\,$ is an open subset of $\,A\,$ with $\,U_i=V_i\cap A_i\,$
for $\,1\leq i\leq n\,$ then $\;U_1\cup\cdots \cup U_n\,=\,
(V_1\cap\cdots\cap V_n)\cup\bigcup_{i=1}^n (V_i\cap(A_i\setminus\{a\})\,$.)
For $\;1\leq i\leq n\;$ consider  $\, A_i\,$ equipped
with a suitable metric $\,d_i\,$. Define a mapping from 
$\,A\times A\,$ into $\,{\Bbb R}\,$ in the following way.
If $\,x,y\in A_i\,$ for some $\,i\,$ then put 
$\,d(x,y)=d_i(x,y)\,$. If $\,x\in A_i\,$ and $\,y\in A_j\,$
for distinct $\,i,j\,$ then put $\,d(x,y)=d_i(x,a)+d_j(y,a)\,$. 
Of course, $\,d\,$ is a metric on the set $\,A\,$.
(One may regard $\,A\,$ as a hedgehog with body $\,a\,$
and spines $\,A_1,...,A_n\,$.) 
By considering the open neighborhoods of the point $\,a\,$ in the 
space $\,A\,$ we conclude that
the topology generated by the metric $\,d\,$ coincides 
with the topology of the space $\,A\,$, {\it q.e.d.}
\medskip
Now we are ready to prove Theorem 2.
Assume $\;{\bf c}\leq\kappa<\lambda\leq 2^\kappa\,$.
(We ignore the case $\,\kappa=\lambda\,$ because this 
case is covered by Proposition 2.)
Let $\,S\,$ be a discrete space of size $\,\kappa\,$
and $\,{\cal F}\,$ an $\omega$-free filter 
on $\,S\,$ with $\,\chi({\cal F})=\lambda\,$. 
Consider the metrizable product space $\;S\times[0,1[\;$ and 
fix $\;p\not\in S\!\times\![0,1[\;$ and
define a topological space $\,\Phi[{\cal F}]\,$ in the following way.
The points in the space $\,\Phi[{\cal F}]\,$ are the elements 
of $\;\{p\}\cup(S\times[0,1[)\;$ and a subset $\,U\,$ of 
$\;\{p\}\cup(S\times[0,1[)\;$ is open if and only if firstly 
$\;U\setminus\{p\}\;$ is open in the product space $\;S\times[0,1[\;$ and 
secondly the point $\,p\,$ lies in 
$\,U\,$ only if 
\smallskip
\centerline{$\;(S\times[t,1[)\cup(F\times[0,1[)\;\subset\; U\;$}
\smallskip
for some $\,t\in[0,1[\;$ and some $\,F\in{\cal F}\,$.
\smallskip
It is plain that this is a correct definition of a topological space
such that the subspace 
$\;\Phi[{\cal F}]\setminus\{p\}\;$ is identical with 
the product space $\,S\times[0,1[\,$. 
Similarly as above we picture $\,\Phi[{\cal F}]\,$ as a cone with apex 
$\,p\,$ and the rulings $\;\{p\}\cup(\{x\}\times[0,1[)\;(x\in X)\;$ 
homeomorphic to the unit interval $\,[0,1]\,$. (Obviously, 
the topology of $\,\Phi[{\cal F}]\,$ is strictly coarser
than the topology of the cone $\,{\cal Q}^*(S)\,$.)
It is straightforward to verify that
$\,\Phi[{\cal F}]\,$
is a regular space. Hence by Corollary~2 the space $\,\Phi[{\cal F}]\,$
is almost metrizable. (Since $\,{\cal F}\,$ is $\omega$-free
and $\,[0,1]\,$ is second countable,
it is clear that $\,\{p\}\,$ is a G$_\delta$-set.)
Since the subspace $\;\{p\}\cup(\{s\}\times[t,1[)\;$
of $\,\Phi[{\cal F}]\,$
is a homeomorphic copy of the compact unit interval $\,[0,1]\,$
for every $\,s\in S\,$ and every $\,t\in[0,1[\,$ 
and since $\,S\,$ is discrete, it is clear that 
$\,\Phi[{\cal F}]\,$ is pathwise connected and locally pathwise connected.
Trivially, $\,|\Phi[{\cal F}]|=\kappa\,$.
\smallskip
Clearly, if $\,{\cal B}\,$ is a filter base on $\,S\,$
generating the filter $\,{\cal F}\,$ then 
\smallskip
\centerline{$\;\big\{\,\{p\}\cup((S\setminus F)\times]1-2^{-n},1[)\cup(F\times[0,1[)
\;\,\big|\,\;n\in{\Bbb N}\,,\;F\in{\cal B}\,\big\}\;$}
\smallskip
is a local basis at $\,p\,$
in the space $\,\Phi[{\cal F}]\,$. 
Conversely, if $\,{\cal U}_p\,$ is a local basis at $\,p\,$
and if we choose for every $\,U\in{\cal U}_p\,$
a real number $\;t_U\in[0,1[\;$ and a set $\;F_U\in{\cal F}\;$
such that $\;(S\times[t_U,1[)\cup(F_U\times[0,1[)\,\subset\,U\;$
then $\;\{\,F_U\;|\;U\in{\cal U}_p\,\}\;$ is a filter base on $\,S\,$
generating the filter $\,{\cal F}\,$. 
Consequently, $\;\chi(p,\Phi[{\cal F}])=\chi({\cal F})\,$.
Therefore, since $\,w(S\times[0,1[)=\kappa\,$, 
we have $\;w(\Phi[{\cal F}])=\chi({\cal F})=\lambda\,$. 
\smallskip
Now consider 
the pathwise connected, locally pathwise connected, 
amost metrizable space $\;\Phi[{\cal F}]\;$
for each of the $\,2^{\lambda}\,$
$\omega$-free filters $\,{\cal F}\,$ on $\,S\,$ with 
$\,\chi({\cal F})=\lambda\,$.
Since the size of each space is $\,\kappa\,$ and 
the weight of each space is $\,\lambda\,$,
by the same arguments about the size of equivalence classes as in the proof of 
Theorem~1 (for $\,\mu=\kappa\,$), the statement in Theorem~2 is true
in case that $\,2^\lambda>2^\kappa\,$
because it is evident that the topologies of the spaces
$\;\Phi[{\cal F}_1]\;$ and $\;\Phi[{\cal F}_2]\;$
are distinct topologies on the set
$\;\{p\}\cup(S\times[0,1[)\;$
whenever $\,{\cal F}_1\,$ and $\,{\cal F}_2\,$ are distinct $\omega$-free 
filters on $\,S\,$.
\medskip
Now assume $\,2^\lambda=2^\kappa\,$ and let 
$\;{\cal P}_\kappa\;$ be a family as provided by 
Proposition 2. Choose one $\omega$-free filter
$\,{\cal F}\,$ on $\,S\,$ with $\,\chi({\cal F})=\lambda\,$
and consider the space $\,\Phi[{\cal F}]\,$. Note that 
$\,x\in\Phi[{\cal F}]\,$ is a noncut point of $\,\Phi[{\cal F}]\,$ 
if and only if $\,x=(s,0)\,$ for some $\,s\in S\,$.
For every $\,H\in{\cal P}_\kappa\,$ create a space $\,X(H)\,$
in the following way. Consider the compact unit square $\,[0,1]^2\,$
and choose a point $\,a_1\in[0,1]^2\,$. (Clearly, $\,a_1\,$ is a
noncut point of $\,[0,1]^2\,$. Note also that no
connected open subset of $\,[0,1]^2\,$ has cut points.)
Choose a noncut point 
$\,a_2\,$ in $\,\Phi[{\cal F}]\,$ and a noncut point 
$\,a_3\,$ in $\,H\,$. Finally, let $\,X(H)\,$ be the 
quotient of the topological sum of the three spaces 
$\,[0,1]^2\,$ and $\,\Phi[{\cal F}]\,$ and $\,H\,$ 
by the subspace $\;\{a_1,a_2,a_3\}\,$.
Roughly speaking, $\,X(H)\,$ is created by 
sticking together the three spaces 
so that the three points $\;a_1,a_2,a_3\;$ are identified.
It is clear that 
$\,X(H)\,$ is pathwise connected and locally pathwise connected
and regular and $\,|X(H)|=\kappa\,$ and $\,w(X(H))=\lambda\,$.
\smallskip
There is precisely one point 
$\,b\in X(H)\,$ with $\,\chi(b,X(H))=\lambda\,$. 
This point $\,b\,$ corresponds with the point $\,p\in \Phi[{\cal F}]\,$.
By virtue of Lemma 4 for $\,n=3\,$ the subspace $\,X(H)\setminus\{b\}\,$
of $\,X(H)\,$ is metrizable.
Consequently, if $\;H\in{\cal P}_\kappa\;$ then
$\,X(H)\,$ is almost metrizable.
The $\,2^\kappa\,$ spaces $\;X(H)\;(H\in{\cal P}_\kappa)\;$
are mutually non-homeomorphic
because each $\,H\in{\cal P}_\kappa\,$ can be recovered from $\,X(H)\,$ 
as follows. 
\smallskip
Since cut points in $\,H\,$ resp.~in $\,\Phi[{\cal F}]\,$ lie dense
and since $\,[0,1]^2\,$ has no cut points, 
there is precisely one point $\,q\,$ in $\,X(H)\,$
such that every neighborhood of $\,q\,$ contains two
nonempty connected open sets $\,U_1,U_2\,$ where
$\,U_1\,$ has no cut points and where $\,U_2\,$ has cut points.
(This point $\,q\,$ must be the point obtained by
identifying the three points $\;a_1,a_2,a_3\,$.)
The subspace $\;X(H)\setminus\{q\}\;$ has precisely three
components and every component of $\;X(H)\setminus\{q\}\;$
is homeomorphic either with $\,\Phi[{\cal F}]\setminus\{a_2\}\,$
or with $\,H\setminus\{a_3\}\,$ or with $\,[0,1]^2\setminus\{a_1\}\,$. 
Therefore, precisely one component is not metrizable.
(If $\,s\in S\,$ then the space 
$\,\Phi[{\cal F}]\setminus\{(s,0)\}\,$ is not metrizable 
since it has no countable local basis at $\,p\,$.)
The two metrizable components of 
$\;X(H)\setminus\{q\}\;$ can be distinguished by 
the observation that one component has infinitely many cut points
while the other component has no cut points.                         
If $\,M\,$ is a metrizable       
component of $\;X(H)\setminus\{q\}\;$ which has cut points 
then the subspace
$\;M\cup\{q\}\;$ of $\,X(H)\,$ is homeomorphic with $\,H\,$, {\it q.e.d.}
\bigskip
{\bff 7. Proof of Theorem 3}
\medskip
{\bf Lemma 5.} {\it There exists a second countable,
countably infinite Hausdorff space
$\,H\,$ such that 
$\,H\setminus E\,$ is connected and locally connected 
for every finite set $\,E\,$.}
\smallskip
{\it Proof.} Let $\,H\,$ be the set $\,{\Bbb N}\,$ equipped with the 
coarsest topology such that if 
$\,p\,$ is a prime and $\,a\in{\Bbb N}\,$ is not divisible by $\,p\,$ then
$\;{\Bbb N}\cap\{\,p+ka\;|\;k\in{\Bbb Z}\,\}\;$ is open.                    
Referring to [13] Nr.~61, $\,H\,$ is a locally connected 
Hausdorff space such that
the intersection of the closures of any two 
nonempty open subsets of $\,H\,$ must be an infinite set.
Therefore, if $\,E\,$ is a finite set then the subspace
$\,H\setminus E\,$ of $\,H\,$ is connected.
Since $\,H\,$ is locally connected,
$\,H\setminus E\,$ is locally connected for every finite set $\,E\,$,
{\it q.e.d.}
\bigskip
The first step in proving Theorem 3 is a proof of the following 
enumeration theorem about countable connected spaces.
\medskip
{\bf Theorem 5.} {\it For every $\;\lambda\leq {\bf c}\;$ 
there exist $\,2^\lambda\,$ mutually non-homeomorphic 
connected, locally connected Hausdorff spaces 
of size $\,\aleph_0\,$ and weight $\,\lambda\,$.}
\medskip
{\it Proof.} Let $\,H\,$ be a connected, locally connected
Hausdorff space with $\;|H|=w(H)=\aleph_0\;$
as provided by Lemma 5.
Fix $\,e\in H\,$ and note that $\,e\,$ is a noncut 
point in $\,H\,$. Put $\;M\,:=\,H\setminus\{e\}\,$.
So $\,M\,$ is connected as well. 
\smallskip
Let $\,S\,$ be an infinite discrete space 
and let $\,{\cal F}\,$ be a free filter on $\,S\,$
with $\,\chi({\cal F})\geq |S|\,$. 
Consider the product space $\;S\times M\;$ and 
fix $\;p\,\not\in\, S\times M\;$ and
consider $\;\Psi[{\cal F}]\,:=\,\{p\}\cup(S\times M)\;$ 
equipped with the following topology.
A subset $\,U\,$ of 
$\;\{p\}\cup(S\times M)\;$ is open if and only if 
$\;U\setminus\{p\}\;$ is open in the product space $\;S\times M\;$ and 
$\,p\in U\,$ implies that 

\centerline{$\;(S\times (V\setminus\{e\}))\cup(F\times M)\;\subset\; U\;$}

for some neighborhood $\,V\,$ of $\,e\,$ in $\,H\,$
and some $\,F\in{\cal F}\,$.
Similarly as in the proof of Theorem 2,
$\;\Psi[{\cal F}]\;$ is a connected and locally
connected Hausdorff space and $\;|\Psi[{\cal F}]|=|S|\;$
and $\;w(\Psi[{\cal F}])=\chi({\cal F})\,$. 
\smallskip
Now let $\,S\,$ be the discrete Euclidean space $\,{\Bbb N}\,$.
If $\;2^\lambda>{\bf c}\;$ then 
with the help of $\,2^\lambda\,$ free filters on $\,{\Bbb N}\,$
with $\,\chi({\cal F})=\lambda\,$ we can track down 
$\,2^\lambda\,$ mutually non-homeomorphic spaces 
$\;\Psi[{\cal F}]\,$. (Note that there are only $\,{\bf c}\,$ permutations 
on $\,{\Bbb N}\,$ and use the argument on sizes of equivalence classes.) 
So it remains to settle the case $\;2^\lambda={\bf c}\,$.
\smallskip
Let $\,Z\,$ be the space $\,\Psi[{\cal F}]\,$ for some
free filter $\,{\cal F}\,$ on $\,{\Bbb N}\,$ 
with $\,\chi({\cal F})=\lambda\,$.
So the underlying set
of $\,Z\,$ is $\;\{p\}\cup ({\Bbb N}\times (H\setminus\{e\}))\;$
and the countable Hausdorff space $\,Z\,$ is connected and locally connected
and $\,w(Z)=\lambda\,$ due to $\,\chi(p,Z)=\lambda\,$.
The point $\,p\,$ is the only cut point of $\,Z\,$
and $\;Z\setminus \{p\}\;$ has infinitely many components.
Keep in mind that $\;|H|=w(H)=\aleph_0\;$ 
and that if $\,a\in H\,$ then the spaces $\,H\,$ and 
$\;H\setminus\{a\}\;$ and $\;H\setminus\{a,e\}\;$
are connected and locally connected. 
Fix  $\;b\,\in\,H\setminus\{e\}\;$
and consider the subset $\;\hat Z\,:=\,\{\,(s,b)\;|\;s\in {\Bbb N}\,\}\;$
of $\,Z\,$.  
Clearly, $\,\hat Z\,$ is closed and discrete and 
$\;Z\setminus\{z\}\;$
is connected and locally connected for every $\,z\in \hat Z\,$.
Choose for every $\,m\in{\Bbb N}\,$ and every $\;i\in\{1,...,m\}\;$
spaces $\;H^{(m)}_i\;$ 
such that $\;H^{(m)}_i\;$ is homeomorphic with $\,H\,$
and $\;H^{(m)}_i\cap H^{(n)}_j=\emptyset\;$
whenever $\,m\not=n\,$ or $\,i\not=j\,$.
Furthermore assume that $\;H^{(m)}_i\cap Z=\emptyset\;$
for every $\,m\,$ and every $\,i\,$.
Let $\,\varphi\,$ be a choice function on the class
of all infinite sets, i.e.~$\;\varphi(A)\in A\;$
for every infinite set $\,A\,$.
Now define for every nonempty set $\,T\subset {\Bbb N}\,$ a Hausdorff space 
$\,Q[T]\,$ as follows. Consider the topological sum $\,\Sigma[T]\,$
of countably infinite and mutually disjoint 
spaces where the summands are 
$\,Z\,$ and all spaces $\;H^{(m)}_i\;$ with $\,m\in T\,$ and 
$\;i\in\{1,...,m\}\,$. 
Define an equivalence relation on 
$\,\Sigma[T]\,$ such that the non-singleton equivalence 
classes are precisely the sets 
$\;\{(m,b)\}\,\cup\,\{\,\varphi(H^{(m)}_1),...,\varphi(H^{(m)}_m)\,\}\;$
with $\,m\in T\,$. 
(Note that $\,(m,b)\in\hat Z\,$ for every $\,m\in T\,$.)
Finally, let $\,Q[T]\,$ denote the quotient 
space of $\,\Sigma[T]\,$ with respect to this equivalence relation.
Roughly speaking, $\,Q[T]\,$ is the union of $\,Z\,$
and all spaces $\;H^{(m)}_i\;$ with $\,m\in T\,$ and $\,i\in\{1,...,m\}\,$
where for every $\,m\in T\,$ the $\,m+1\,$ points 
$\;(m,b),\varphi(H^{(m)}_1),...,\varphi(H^{(m)}_m)\;$
are identified. We consider $\,Z\,$ to be a subset of $\,Q[T]\,$.
One may picture $\,Q[T]\,$ as an expansion of $\,Z\,$
created by attaching $\,m\,$ copies
of $\,H\,$ to $\,Z\,$ 
at the point $\,(m,b)\in \hat Z\,$ for every $\,m\in T\,$.
It is evident that 
$\,Q[T]\,$ is a connected and locally connected countably infinite 
Hausdorff space. We have $\;w(Q[T])=\lambda\;$
since $\,Z\,$ is a subspace of $\,Q[T]\,$ 
with $\,w(Z)=\lambda\,$ 
and $\;\chi(x, Q[T])=\aleph_0\;$ if $\;p\not=x\in Q[T]\,$.
Thus the case $\,2^\lambda={\bf c}\,$ is settled by verifying that 
two spaces $\,Q[T_1]\,$ and $\,Q[T_2]\,$ cannot be homeomorphic 
if $\;\emptyset\not=T_1,T_2\subset{\Bbb N}\;$ and $\,T_1\not= T_2\,$.
This must be true because the set $\,T\subset{\Bbb N}\,$
is completely determined by the topology of $\,Q[T]\,$
by the following observation.
\smallskip
Let $\,\emptyset\not=T\subset{\Bbb N}\,$. For every point $\,x\in Q[T]\,$ let 
$\,\nu(x)\,$ denote the total number of all components of 
the subspace $\;Q[T]\setminus\{x\}\,$. 
The following three statements for $\,x\in Q[T]\,$ are evident.
Firstly, $\;\nu(x)\geq \aleph_0\;$
if and only if $\,x=p\,$.
Secondly, $\;1<\nu(x)<\aleph_0\;$
if and only if $\;x=(m,b)\in\hat Z\;$ for some $\,m\in T\,$.
Thirdly, $\;\nu(x)=1\;$
if and only if $\,x\,$ is an element of the set 
$\;Q[T]\setminus((T\times\{b\})\cup\{p\})\,$.
Concerning the second statement we compute 
$\;\nu((m,b))=m+1\;$ for every  
$\,m\in T\,$. Consequently, 
\smallskip
\centerline{$\;\{\,\nu(x)-1\;\,|\,\;x\in Q[T]\;\land\;\nu(x)\in {\Bbb N}\,\}
\setminus\{0\}\;=\;T\;$}
\smallskip
whenever $\,T\,$ is one of the $\,{\bf c}\,$ non-empty subsets of $\,{\Bbb N}\,$, 
{\it q.e.d.}
\bigskip
Now in order to prove Theorem 3 assume $\;\aleph_0\leq\kappa<{\bf c}\;$
and $\;\kappa\leq\lambda\leq 2^\kappa\,$. Referring to Theorem 5
there is nothing more to show in case that $\,\kappa=\aleph_0\,$.
So we also assume that $\,\kappa>\aleph_0\,$.
Let $\,S\,$ be a discrete space of size $\,\kappa\,$.
By Proposition 3 there are $\,2^\lambda\,$ 
free filters $\,{\cal F}\,$ on $\,S\,$
with $\,\chi({\cal F})=\lambda\,$. 
For each one of these filters $\,{\cal F}\,$ 
consider the connected and locally connected Hausdorff 
space $\,\Psi[{\cal F}]\,$ of size $\,\kappa\,$ and weight $\,\lambda\,$
as defined in the previous proof. 
Hence in case that 
$\,2^\lambda>2^\kappa\,$ we can track down $\,2^\lambda\,$ 
filters $\,{\cal F}\,$ on $\,S\,$ such that the corresponding spaces
$\,\Psi[{\cal F}]\,$ are mutually non-homeomorphic.
\medskip
So it remains to settle the case
$\;2^\lambda=2^\kappa\,$. 
Choose any free filter $\,{\cal F}\,$ on $\,S\,$ with
$\,\chi({\cal F})=\lambda\,$ and
and consider the space $\;\Psi:=\Psi[{\cal F}]\;$ 
of size $\,\kappa\,$ and weight
$\,\lambda\,$. Fix a noncut point
$\,z\in \Psi\,$. Keep in mind that 
$\,\Psi\,$ has precisely one cut point $\,p\,$
and that $\;\Psi[{\cal F}]\setminus \{p\}\;$ has precisely $\,\kappa\,$ 
and hence {\it uncountably many} components. 
\medskip
In view of our proof of Theorem 5 there is a family $\,{\cal C}\,$ 
of mutually non-homeomorphic {\it countable} Hausdorff 
spaces of weight $\,\kappa\,$
(and hence not necessarily of weight $\,\lambda\,$)
such that $\,|{\cal C}|=2^\kappa\,$  
and if $\,C\in{\cal C}\,$ then $\,C\,$ is connected and locally connected
and contains precisely one cut point $\,q(C)\,$ such that 
$\;C\setminus \{q(C)\}\;$ has infinitely many components.  
In particular, all these components are countable and $\,\aleph_0\,$ is
their total number.
\medskip
For every $\,C\in {\cal C}\,$ consider the topological sum 
$\,\Psi+C\,$ and define an equivalence relation 
such that $\,\{z,q(C)\}\,$ is an 
equivalence class and all other equivalence classes are singletons.
Let $\,Q[C]\,$ denote the quotient 
space of $\;\Psi+C\;$ with respect to this equivalence relation.
So $\,Q[C]\,$ is obtained by sticking together the 
spaces $\,\Psi\,$ and $\,C\,$ at one point 
and this point is the identification of $\,z\in \Psi\,$
and $\,q(C)\in C\,$. It is clear that $\,Q[C]\,$ 
is a connected, locally connected 
Hausdorff space of size $\,\kappa\,$ and weight
$\,\lambda\,$. So we are done by verifying 
that for distinct $\,C_1,C_2\in{\cal C}\,$
the spaces $\,Q[C_1]\,$ and $\,Q[C_2]\,$ are never homeomorphic.
This must be true because each $\,C\in{\cal C}\,$ can be 
recovered from $\,Q[C]\,$ as follows. 
\medskip
There is a unique point $\,\xi\,$ in $\,Q[C]\,$
such that $\;Q[C]\setminus \{\xi\}\;$ has precisely 
$\,\aleph_0\,$ components. (This point $\,\xi\,$ is the one corresponding
with the equivalence class $\,\{z,q(C)\}\,$.)
Among these components there is precisely one of uncountable size.
(This component is the one which contains the point $\,p\in\Psi\,$.)
Let $\,K\,$ be the unique uncountable
component of $\;Q[C]\setminus \{\xi\}\,$.
Then $\;Q[C]\setminus K\;$ 
is essentially identical, at least homeomorphic with the space $\,C\,$.
\bigskip
{\bff 8. Proof of Theorem 4}
\medskip
Our goal is to derive Theorem 4 from Theorem 1
by using appropriate modifications of the cones $\,{\cal Q}^*(X)\,$
considered in Section 6. In order to accomplish this we need
building blocks provided by the following lemma.
\medskip
{\bf Lemma 6.} {\it There exists a second countable, connected,
totally pathwise disconnected, nowhere locally connected, 
metrizable space $\,M\,$ of size $\,{\bf c}\,$ 
which contains precisely one noncut point $\,b\,$ and where 
$\,M\setminus\{x,b\}\,$ has precisely two components whenever
$\,b\not=x\in M\,$.}
\medskip
{\it Proof.} Let $\,f\,$ be a function 
from $\,{\Bbb R}\,$ into $\,{\Bbb R}\,$ such that 
the graph of $\,f\,$ is a 
dense and connected subset of the Euclidean plane $\,{\Bbb R}^2\,$.
(See [10] for a construction of such a function $\,f\,$.) Automatically,
$\,f\,$ is discontinuous everywhere.
Let $\,M\,$ be the intersection of $\;[0,\infty[\times{\Bbb R}\;$
and the graph of $\,f\,$. It is straightforward to check 
that $\,M\,$ fits with $\,b=(0,f(0))\,$, {\it q.e.d.}
\medskip
Now we are ready to prove Theorem 4.
Assume $\;{\bf c}\leq \kappa\leq \lambda\leq 2^\kappa\;$ and
let $\,{\cal Y}={\cal Y}(\kappa,\lambda)\,$ be a family of precisely 
$\,2^\lambda\,$ mutually non-homeomorphic scattered, normal
spaces of size $\,\kappa\,$ and weight $\,\lambda\,$
such that if $\,Y\in{\cal Y}\,$ then
for a certain finite set $\,\gamma(Y)\subset Y\,$ 
the subspace $\;Y\setminus\gamma(Y)\;$ is metrizable 
(and hence of weight $\,\kappa\,$) and 
$\,\gamma(Y)\,$ is a G$_\delta$-set in $\,Y\,$.
Precisely, the set $\,\gamma(Y)\,$ is empty when 
$\,\kappa=\lambda\,$ and a singleton $\,\{y\}\,$ 
when $\,\kappa<\lambda\,$. (Clearly, if $\,\gamma(Y)=\{y\}\,$
then $\,\chi(y,Y)=\lambda\,$.)
If $\,2^\lambda>2^\kappa\,$ then such a family $\,{\cal Y}\,$ 
exists by considering the $\,2^\lambda\,$ almost discrete spaces
provided by Theorem 1. 
If $\,\lambda>\kappa\,$ and $\,2^\lambda=2^\kappa\,$ then such a
family $\,{\cal Y}\,$ exists in view of the construction 
in Section 5 which proves Theorem 1 in case that $\,2^\lambda=2^\kappa\,$.
If $\,\lambda=\kappa\,$ then such a family $\,{\cal Y}\,$ 
exists by Proposition 1.
\smallskip
Let $\,M\,$ be a metrizable space as in Lemma 6 
and let $\,b\,$ denote the noncut point of $\,M\,$
and fix a point $\;a\,\in\,M\setminus\{b\}\,$.
For an infinite, scattered, normal space $\,X\,$ consider the product space 
$\;X\times M\;$ and 
fix $\;p\not\in X\!\times\!M\;$ and
put $\;K(X)\,:=\,\{p\}\cup(X\times(M\setminus\{b\})\,$. 
Declare a subset $\,U\,$ of $\,K(X)\,$ open if and only if 
$\;U\setminus\{p\}\;$ is open in the product space 
$\;X\times(M\setminus\{b\})\;$ and 
$\,p\in U\,$ implies that $\,U\,$ contains
$\;X\times(N\setminus\{b\})\;$ for some 
neighborhood $\,N\,$ of $\,b\,$ in the space $\,M\,$.
It is plain that $\,K(X)\,$ is a well-defined regular space.
Since $\,M\,$ is metrizable and 
$\,\chi(p,K(X))=\aleph_0\,$, if $\,X\,$ is metrizable then
$\,K(X)\,$ has a $\sigma$-locally finite base and hence
$\,K(X)\,$ is metrizable.
\smallskip
Now for $\,Y\in{\cal Y}\,$ consider the subspace 
$\;L(Y)\;:=\;K(Y)\setminus(\gamma(Y)\times(M\setminus\{a,b\}))\;$
of $\,K(Y)\,$ and the subspace 
$\;S(Y)\;:=\;L(Y)\setminus(\gamma(Y)\times\{a\})\;$
of $\,L(Y)\,$. Trivially, the spaces $\,K(Y)\,$ and $\,L(Y)\,$
and $\,S(Y)\,$ coincide if $\,\kappa=\lambda\,$.
Furthermore the space $\,S(Y)\,$ 
coincides with the metrizable space $\;K(Y\setminus\gamma(Y))\,$.
Therefore and by Corollary 2, $\,L(Y)\,$ is an almost metrizable space 
since $\;\gamma(Y)\times\{a\}\;$ is a G$_\delta$-set in $\,K(Y)\,$
of size $\,0\,$ or $\,1\,$.
We have $\,|L(Y)|=\kappa\,$ and $\,w(L(Y))=\lambda\,$
because if $\,\gamma(Y)=\{y\}\,$ then $\;\chi((y,a),L(Y))=\lambda\,$.
It is evident that $\,S(Y)\,$ is 
connected and totally pathwise disconnected and nowhere locally connected.
Consequently, $\,L(Y)\,$ 
is totally pathwise disconnected and nowhere locally connected.
And $\,L(Y)\,$ is connected since the connected set $\,S(Y)\,$ is 
dense in $\,L(Y)\,$.
\smallskip
Finally, the spaces $\;L(Y)\;(Y\in{\cal Y})\;$ are mutually non-homeomorphic
because every $\,Y\in{\cal Y}\,$ can be recovered from $\,L(Y)\,$.
Indeed, for $\,x\in L(Y)\,$ let $\,{\cal C}(x)\,$ denote the
family of all components of the subspace $\;L(Y)\setminus\{x\}\;$ of 
$\,L(Y)\,$. Then $\,{\cal C}(x)\,$ is an infinite set 
if and only if $\,x=p\,$.
Because the scattered space $\,Y\,$ has infinitely many 
isolated points and if $\,u\in Y\,$ is isolated then 
$\;\{u\}\times(M\setminus\{b\})\;$ lies in $\,{\cal C}(x)\,$.
If $\,u\,\in Y\setminus\gamma(Y)\,$ and $\,b\not=v\in M\,$
then $\,|{\cal C}((u,v))|\leq 2\,$ (and
$\,|{\cal C}((u,v))|=2\,$ when $\,u\,$ is isolated in $\,Y\,$).
And if $\,\gamma(Y)=\{y\}\,$ then $\,|{\cal C}((y,a))|=1\,$.
Thus $\;\{p\}\,=\,\{\,x\in L(Y)\;|\;|{\cal C}(x)|\geq\aleph_0\,\}\,$,
whence the point $\,p\,$ can be recovered from the space $\,L(Y)\,$.
Now let $\,{\cal C}\,$ be the family of all components of the space 
$\;L(Y)\setminus\{p\}\,$. Since $\,Y\,$ is totally disconnected, 
the members of $\,{\cal C}\,$ are precisely the sets 
$\;\{u\}\times(M\setminus\{b\})\;$ with $\;u\,\in\,Y\setminus\gamma(Y)\;$
plus the singleton $\,\gamma(Y)\times\{a\}\,$ if and only if 
$\,\gamma(Y)\not=\emptyset\,$. Naturally, the quotient space 
of $\;L(Y)\setminus\{p\}\;$ by the equivalence relation defined 
via the partition $\,{\cal C}\,$ is homeomorphic with $\,Y\,$
for every $\,Y\in{\cal Y}\,$.
This concludes the proof of Theorem 4.
\bigskip
{\bff 9. Overweight compact spaces}
\medskip
While $\,w(X)\leq |X|\,$
for every compact Hausdorff space $\,X\,$ (see [2, 3.3.6]),
for compact T$_1$-spaces $\,X\,$ one cannot rule out 
$\,w(X)>|X|\,$ and actually we can prove the following enumeration theorem
by applying Theorems 1 and 2 and 3. 
\medskip
{\bf Theorem 6.} {\it If $\;\kappa\leq\lambda\leq 2^\kappa\;$
then there exist two families $\,{\cal C}_1,{\cal C}_2\,$
of mutually non-homeomorphic compact 
{\rm T}$_1$-spaces of size $\,\kappa\,$ and weight $\,\lambda\,$
such that $\;|{\cal C}_1|=|{\cal C}_2|=2^\lambda\;$
and all spaces in $\,{\cal C}_1\,$ are scattered,
all spaces in $\,{\cal C}_2\,$ are connected 
and locally connected, and if $\,\kappa\geq{\bf c}\,$ 
then all spaces in $\,{\cal C}_2\,$ are arcwise connected 
and locally arcwise connected.}
\medskip
In order to prove Theorem 6 we consider T$_1$-compactifications 
of Hausdorff spaces. If $\,Y\,$ is an infinite Hausdorff space 
with $\;|Y|\leq w(Y)\;$
then define a topological space $\,\Gamma(Y)\,$ which expands $\,Y\,$
in the following way. Put $\;\Gamma(Y)\,=\,Y\cup\{z\}\;$
where $\,z\not\in Y\,$ and declare $\,U\subset \Gamma(Y)\,$
open either when $\,U\,$ is an open subset of $\,Y\,$ 
or when $\,z\in U\,$ and $\,Y\setminus U\,$ is finite. 
It is clear that in this way a topology on $\,\Gamma(Y)\,$ is well-defined
such that $\,Y\,$ is a dense subspace of $\,\Gamma(Y)\,$.
Obviously, $\;\Gamma(Y)\setminus\{x\}\;$ 
is open for every $\,x\in \Gamma(Y)\,$
and hence $\,\Gamma(Y)\,$ is a T$_1$-space.
Since all neighborhoods of $\,z\,$ cover 
the whole space $\,\Gamma(Y)\,$ except finitely many points, 
$\,\Gamma(Y)\,$ is compact.
Trivially, $\,|\Gamma(Y)|=|Y|\,$.  
We have $\,w(\Gamma(Y))=w(Y)\,$ since 
$\,w(Y)\geq |Y|\,$ and $\,Y\,$ is a subspace 
of $\,\Gamma(Y)\,$ and, by definition, there is a local basis
at $\,z\,$ of size $\,|Y|\,$.
\smallskip
Evidently, if $\,Y\,$ is scattered then $\,\Gamma(Y)\,$ is scattered.
On the other hand it is clear that if $\,Y\,$ is dense in itself
then $\,\Gamma(Y)\,$ is connected and every neighborhood of $\,z\,$ is 
connected. So if $\,Y\,$ is connected and locally connected then 
$\,\Gamma(Y)\,$ is connected and locally connected. 
\smallskip
We claim that if $\,Y\,$ is pathwise connected 
then $\,\Gamma(Y)\,$ 
is arcwise connected. Assume that the Hausdorff space 
$\,Y\,$ is pathwise connected and hence
arcwise connected and let $\,a\in Y\,$.
Of course it is enough to find an arc which 
connects the point $\,a\,$ with the point $\,z\not\in Y\,$.
Since $\,Y\,$ is arcwise connected,
we can define a homeomorphism $\,\varphi\,$ from $\,[0,1]\,$ onto a subspace 
of $\,Y\,$ such that $\,\varphi(0)=a\,$. Define an injective function $\,f\,$ 
from $\,[0,1]\,$ into $\,\Gamma(Y)\,$ via $\,f(1)=z\,$ and
$\,f(t)=\varphi(t)\,$ for $\,t<1\,$. Let $\,U\,$ be an open subset 
of $\,\Gamma(Y)\,$. If $\,z\in U\,$ then 
$\;U\setminus Y\;$ is finite and thus  
$\;f^{-1}(U)\;$ is a cofinite and hence open subset of $\,[0,1]\,$.
If $\,z\not\in U\,$ then $\,U\,$ is an open subset of $\,Y\,$
and hence $\;f^{-1}(U)=\varphi^{-1}(U)\setminus\{1\}\;$
is an open subset of $\,[0,1]\,$. Thus the injective function 
$\,f\,$ is continuous. 
\smallskip
Since every neighborhood of $\,z\,$ contains all but finitely many 
points from $\,Y\,$, by exactly the same arguments we conclude
that if $\,Y\,$ is locally pathwise connected then 
every neighborhood of $\,z\,$ is an arcwise connected subspace of 
$\,\Gamma(Y)\,$. Consequently, if the Hausdorff space
$\,Y\,$ is locally pathwise connected 
then the T$_1$-space $\,\Gamma(Y)\,$ is locally arcwise connected.
\smallskip
The space $\,Y\,$ can be recovered
from $\,\Gamma(Y)\,$ (up to homeomorphism)
provided that $\,Y\,$ has at least two limit points.
Because then it is evident that 
$\,z\,$ is the unique point
$\,x\in\Gamma(Y)\,$ such that 
the subspace $\;\Gamma(Y)\setminus\{x\}\;$ of $\,\Gamma(Y)\,$ is Hausdorff.
\medskip
By virtue of Theorem 1, for $\,\kappa\leq\lambda\leq 2^\kappa\,$ 
let $\,{\cal Y}_1(\kappa,\lambda)\,$ be a family 
of mutually non-homeomorphic, scattered Hausdorff spaces 
of size $\,\kappa\,$ and weight $\,\lambda\,$ such that 
$\,|{\cal Y}_1(\kappa,\lambda)|=2^\lambda\,$.
By virtue of Theorem 3, for $\,\kappa<{\bf c}\,$ and 
$\,\kappa\leq\lambda\leq 2^\kappa\,$ 
let $\,{\cal Y}_2(\kappa,\lambda)\,$ be a family 
of mutually non-homeomorphic connected and locally connected 
Hausdorff spaces 
of size $\,\kappa\,$ and weight $\,\lambda\,$ such that 
$\,|{\cal Y}_2(\kappa,\lambda)|=2^\lambda\,$.
By virtue of Theorem 2, for $\,{\bf c}\leq\kappa\leq\lambda\leq 2^\kappa\,$ 
let $\,{\cal Y}_3(\kappa,\lambda)\,$ be a family 
of mutually non-homeomorphic pathwise connected and locally 
pathwise connected Hausdorff spaces 
of size $\,\kappa\,$ and weight $\,\lambda\,$ such that 
$\,|{\cal Y}_3(\kappa,\lambda)|=2^\lambda\,$.
Now put
$\;{\cal C}_1\;:=\;\{\,\Gamma(Y)\;|\;Y\in{\cal Y}_1(\kappa,\lambda)\,\}\;$
and $\;{\cal C}_2\;:=
\;\{\,\Gamma(Y)\;|\;Y\in{\cal Y}_i(\kappa,\lambda)\,\}\;$ 
where $\,i=2\,$ when $\,\kappa<{\bf c}\,$ and $\,i=3\,$ when 
$\,\kappa\geq{\bf c}\,$. Then $\,{\cal C}_1,{\cal C}_2\,$ 
are families which prove Theorem 6.
\bigskip
The condition $\,\kappa\geq{\bf c}\,$ in Theorem 6 is inevitable
since, trivially, $\;|X|\geq{\bf c}\;$
for every infinite, {\it arcwise connected} space.
There arises the question whether $\;|X|\geq{\bf c}\;$
is inevitable for infinite, {\it pathwise connected} T$_1$-spaces.
(Of course, every finite T$_1$-space $\,X\,$
is discrete and hence not pathwise connected when $\,|X|\geq 2\,$.)
It is well-known that a pathwise connected T$_1$-space of size 
$\,\aleph_0\,$ does not exist (see also Proposition 5 below).
So the essential question is whether
there are pathwise connected T$_1$-spaces $\,X\,$ with 
$\,\aleph_0<|X|<{\bf c}\,$ (provided that there are cardinals $\,\mu\,$
with $\,\aleph_0<\mu<{\bf c}\,$). The following 
proposition shows that there is no chance to track down 
such spaces $\,X\,$.
\medskip
{\bf Proposition 5.} {\it Pathwise connected
{\rm T}$_1$-spaces $\,X\,$ with $\,2\leq |X|\leq \aleph_0\,$ do not exist.}
{\it It is consistent with {\rm ZFC} that 
$\;|\{\,\kappa\;|\;\aleph_0<\kappa<{\bf c}\,\}|>\aleph_0\;$
and pathwise connected
{\rm T}$_1$-spaces $\,X\,$ with $\,\aleph_0<|X|<{\bf c}\,$ do not exist.}
\medskip
If $\,X\,$ is a T$_1$-space and 
$\;f\,:[0,1]\to X\;$ is continuous then 
$\;\{\,f^{-1}(\{x\})\;|\;x\in X\,\}\setminus\{\emptyset\}\;$
is a decomposition of $\,[0,1]\,$ into precisely $\,|f([0,1])|\,$ 
nonempty closed subsets. Therefore, Proposition 5 is an immediate 
consequence of 
\medskip

{\bf Proposition 6.} {\it Every partition of $\,[0,1]\,$ into 
at least two closed sets is uncountable. It is consistent 
with {\rm ZFC} that 
uncountably many cardinals 
$\,\kappa\,$ with $\,\aleph_0<\kappa<{\bf c}\,$ exist 
while still a partition $\,{\cal P}\,$ of $\,[0,1]\,$ 
into closed sets with $\,\aleph_0<|{\cal P}|<{\bf c}\,$ does not exist.}
\medskip
Certainly, the first statement in Proposition 6 is an 
immediate consequence of Sierpi\'nski's theorem [2, 6.1.27]. 
However, in order to prove Proposition 6 we need another approach than 
in the proof of [2, 6.1.27]. (Moreover, the following proof is  
much easier than the proof of [2, 6.1.27].)
\medskip
Assume that $\,{\cal P}\,$ is a partition 
of $\,[0,1]\,$ into closed sets with $\,|{\cal P}|\geq 2\,$.
For $\,S\subset [0,1]\,$ let $\,\partial S\,$ 
denote the boundary of $\,S\,$ in the compact space $\,[0,1]\,$. 
(Notice that then $\,\partial[0,1]=\emptyset\,$.)
Put $\;{\cal V}\,:=\,\{\,\partial A\;|\;A\in{\cal P}\,\}\;$
and $\;W\,:=\,\bigcup{\cal V}\,$.
Then $\,\emptyset\not\in{\cal V}\,$ since $\,[0,1]\not\in{\cal P}\,$
and hence $\,{\cal V}\,$ is a partition of $\,W\,$ with
$\,|{\cal V}|=|{\cal P}|\,$.
The nonempty set $\,W\,$ is a closed subset of $\,[0,1]\,$ because 
$\;W\,=\,[0,1]\setminus\bigcup\{\,A\setminus\partial A\;|\;A\in{\cal P}\,\}\;$
since $\,{\cal P}\,$ is a partition of $\,[0,1]\,$.
We claim that the closed sets $\,V\in{\cal V}\,$  
are nowhere dense in the compact metrizable space $\,W\,$.
\smallskip
Let $\,A\in{\cal P}\,$ and
assume indirectly that $\,a\,$ is an interior 
point of $\,\partial A\,$ in $\,W\,$.
Then there is an interval $\,I\,$ open in the compact space $\,[0,1]\,$
with $\,a\in I\,$ and $\;I\cap W\subset\partial A\,$. 
Since $\,a\,$ lies in the boundary of $\,A\,$, 
the interval $\,I\,$ intersects $\;[0,1]\setminus A\;$ 
and hence for some $\,B\not=A\,$ in the family $\,{\cal P}\,$
we have $\;I\cap B\not=\emptyset\,$.
However, $\;I\cap \partial B=\emptyset\;$ 
in view of $\,(\partial A)\cap(\partial B)=\emptyset\,$
and $\;I\cap W\subset\partial A\,$. 
Therefore, $\,I\cap B\,$ is a nonempty set 
which is open and closed in the connected space 
$\,I\,$ and hence $\,I\cap B=I\,$ contrarily 
with $\;A\cap I\not=\emptyset\;$
and $\,A\cap B=\emptyset\,$.
\smallskip
Thus $\,{\cal V}\,$ is a partition of the compact Hausdorff space $\,W\,$ 
into nowhere dense subsets with $\,|{\cal V}|=|{\cal P}|\,$.
Therefore $\,|{\cal P}|\leq\aleph_0\,$ is impossible
since $\,W\,$ is a space of second category.
This concludes the proof of the first statement.
Under the assumption of {\it Martin's Axiom} (see [4, 16.11])
also the weaker inequality $\,|{\cal V}|=|{\cal P}|<{\bf c}\,$
is impossible because it is well-known that Martin's axiom
implies that no separable, compact Hausdorff space can be covered by 
less than $\,{\bf c}\,$ nowhere dense subsets.
(Actually, Martin's axiom is {\it equivalent} to
the statement that in every compact Hausdorff space of countable cellularity
any intersection of less than $\,{\bf c}\,$ dense, open sets is dense.)
Therefore, the proof of Proposition 6
is concluded by checking that the existence of 
uncountably many infinite cardinals below $\,{\bf c}\,$
is consistent with ZFC plus Martin's Axiom.
This is certainly true because by applying 
the Solovay-Tennenbaum theorem [4, 16.13]
there is a model of ZFC in which Martin's Axiom holds and
the identity $\,2^{\aleph_0}=\aleph_{\omega_1+1}\,$ is enforced.
(If $\,{\bf c}=\aleph_{\omega_1+1}\,$ then
$\;|\{\,\kappa\;|\;\kappa<{\bf c}\,\}|=\aleph_1>\aleph_0\,$.)
\bigskip
{\it Remark.} There is an interesting observation
concerning compactness and the Hausdorff separation axiom.
By applying Theorem 6 and (1.2), {\it there exist precisely $\,{\bf c}\,$
compact, countable, second countable {\rm T}$_1$-spaces  
up to homeomorphism.} If in this statement T$_1$ is sharpened to 
T$_2$ then we obtain an {\it unprovable hypothesis.}
Indeed, due to Mazurkiewicz and Sierpi\'nski [12], there exist precisely 
$\,\aleph_1\,$ countable (and hence second countable) compact Hausdorff
spaces up to homeomorphism. (The hypothesis 
$\,\aleph_1<{\bf c}\,$ is irrefutable since it is a trivial 
consequence of (1.1).) This discrepancy of provability vanishes 
when {\it uncountable} compacta are counted up to homeomorphism.
Indeed,  by virtue of [7, Theorem 3] 
it can be accomplished that in Theorem 6 for $\,\kappa=\lambda>\aleph_0\,$ 
all spaces in the family $\,{\cal C}_1\,$ 
are Hausdorff spaces.
(Note that $\,w(X)=|X|\,$ for every scattered, compact Hausdorff space.) 
\vfill\eject
\bigskip
{\bff 10. Pathwise connected, scattered spaces}
\medskip
Naturally, a scattered T$_1$-space is totally disconnected and 
hence far from 
being pathwise connected. Furthermore it is plain that no scattered space
is {\it arcwise} connected.
Therefore and in view of Proposition 5
the following enumeration theorem is worth mentioning.
\medskip
{\bf Theorem 7.} {\it If $\;\kappa\leq\lambda\leq 2^\kappa\;$
then there exist two families $\,{\cal C},{\cal L}\,$
of mutually non-homeomorphic 
pathwise connected, scattered {\rm T}$_0$-spaces  
of size $\,\kappa\,$ and weight $\,\lambda\,$
such that $\;|{\cal C}|=|{\cal L}|=2^\lambda\;$
and all spaces in $\,{\cal C}\,$ are compact 
and if $\,\kappa\leq{\bf c}\,$ or $\,2^\kappa<2^\lambda\,$
then all spaces in $\,{\cal L}\,$ are locally pathwise connected.}
\medskip
The existence of the family $\,{\cal C}\,$ in Theorem 7 
can be derived from Theorem 1 in view of the following 
considerations. Let $\,X\,$ be an infinite Hausdorff space.
Fix $\,b\not\in X\,$
and define a topology on the set $\;B(X)\,=\,X\cup\{b\}\;$
by declaring $\;U\subset B(X)\;$ open when either $\,U=B(X)\,$
or $\,U\,$ is an open subset of $\,X\,$. 
(Then $\,\{b\}\,$ is closed and $\,B(X)\,$ is the only neighborhood 
of $\,b\,$.) Obviously, $\,B(X)\,$ is a compact T$_0$-space
and $\,b\,$ is a limit point of 
every nonempty subset of $\;X\,=\,B(X)\setminus\{b\}\,$.
It is trivial that $\,|B(X)|=|X|\,$ and clear that $\,w(B(X))=w(X)\,$.
For any pair $\,x,y\,$ of distinct points in $\,B(X)\,$
define a function $\,f\,$
from $\,[0,1]\,$ into $\,B(X)\,$ via $\,f(t)=x\,$
when $\,t<{1\over 2}\,$ and $\,f({1\over 2})=b\,$ and $\,f(t)=y\,$
when $\,t>{1\over 2}\,$. It is plain that $\,f\,$ is continuous,
whence $\,B(X)\,$ is pathwise connected.
Obviously, if $\,X\,$ is scattered then $\,B(X)\,$ is scattered.
Finally, the space $\,X\,$ can be recovered from $\,B(X)\,$ 
since a singleton $\,\{a\}\,$ is closed in $\,B(X)\,$ 
if and only if $\,a=b\,$. 
\medskip 
Unfortunately, if $\,X\,$ is scattered and not discrete
then $\,B(X)\,$ is not locally connected. Fortunately,
finishing the proof of Theorem 7
we can track down a family $\,{\cal L}\,$ as desired by
adopting the proofs of Theorem 3 and Theorem 5 in Section 7 line by line
such that, throughout, the building block
$\,H\,$ in the definition of $\,\Phi[{\cal F}]\,$ provided by Lemma 5 
is replaced with the space $\,G\,$ provided by the following lemma. 
(In Section 7 the restriction $\,\kappa<{\bf c}\,$  is only for 
avoiding an overlap between Theorem 2 and Theorem 3 and can 
clearly be expanded to $\,\kappa\leq{\bf c}\,$.
The case $\,2^\kappa<2^\lambda\,$ is settled by 
the $\,2^\lambda\,$ spaces $\,\Psi[{\cal F}]\,$ of {\it arbitrary} size 
$\,\kappa\,$.)
\medskip
{\bf Lemma 7.} {\it There exists a second countable, scattered,
countably infinite {\rm T}$_0$-space
$\,G\,$ such that                  
$\,G\setminus E\,$ is pathwise connected and locally pathwise connected
for every finite set $\,E\,$.}
\smallskip
{\it Proof.} Let $\,G\,$ be 
the set $\,\{\,n\in{\Bbb Z}\;|\;n\geq 2\,\}\,$ 
equipped with {\it divisor topology} as defined 
in [13, {\bf 57}]. (A basis of the divisor topology 
is the family of all sets 
$\;\{\,m\in{\Bbb Z}\;|\;m\geq 2\;\land\;m|n\,\}\;$ 
with $\,n\in G\,$.) In view of the considerations in [13], 
it is straightforward to verify that $\,G\,$ fits, {\it q.e.d.}
\medskip\smallskip
{\it Remark.} If $\,i\in\{0,1,2\}\,$
and $\,{\cal F}_i\,$ is a family of mutually non-homeomorphic compact
{\rm T}$_i$-spaces $\,X\,$ with $\,w(X)\leq\kappa\,$ then 
$\,|{\cal F}_i|\leq 2^\kappa\,$ is true for $\,i=2\,$. 
(Because any compact Hausdorff space of weight at most $\,\kappa\,$
is embeddable into the Hilbert cube $\,[0,1]^{\kappa}\,$  
and, since $\,w([0,1]^\kappa)=\kappa\,$ and $\,|X|=2^\kappa\,$,
the compact Hausdorff space 
$\,[0,1]^{\kappa}\,$ has precisely $\,2^\kappa\,$ closed subspaces.)
However, the estimate $\,|{\cal F}_i|\leq 2^\kappa\,$ is false 
for $\,i=0\,$ because $\,|{\cal F}_0|=2^{2^\kappa}\,$ can be achieved
for every $\,\kappa\,$. (In view of (1.2) and since
$\;\max\{|X|,w(X)\}\leq\min\{2^{|X|},2^{w(X)}\}\;$
for every infinite T$_0$-space $\,X\,$, 
$\,2^{2^\kappa}\,$ is the  maximal possible cardinality.)
Indeed, consider for $\,X=[0,1]^\kappa\,$ the compact T$_0$-space 
$\,B(X)=X\cup\{b\}\,$ of size $\,2^\kappa\,$ and weight 
$\,\kappa\,$ defined as above. 
Clearly, for every nonempty $\,S\subset X\,$
the subspace $\,S\cup\{b\}\,$ of $\,B(X)\,$ is compact. 
Since $\,X\,$ is Hausdorff and $\,w(X)=\kappa\,$, there are 
$\,2^{|X|}=2^{2^\kappa}\,$ mutually 
non-homeomorphic subspaces of $\,X\,$
and hence $\,2^{2^\kappa}\,$ mutually 
non-homeomorphic compact subspaces of $\,B(X)\,$.
There arises the interesting 
question whether the estimate $\,|{\cal F}_i|\leq 2^\kappa\,$ 
is generally true for $\,i=1\,$.
\bigskip\smallskip
{\bff 11. Counting P-spaces}
\medskip
A natural modification of the proof of Theorem 1 leads 
to a noteworthy enumeration theorem about $P$-spaces.
As usual (see [1]), a Hausdorff space is a {\it P-space}
if and only if any intersection of countably many open sets is open.
More generally, a Hausdorff space $\,X\,$
is a $P_\alpha$-{\it space} if and only if $\,\alpha\,$
is an infinite cardinal number and $\,\bigcap{\cal U}\,$
is an open subset of $\,X\,$ whenever $\,{\cal U}\,$ 
is a family of open subsets of $\,X\,$ with $\,0\not=|{\cal U}|<\alpha\,$.
So if $\,\alpha=\aleph_0\,$ then every Hausdorff space is a $P_\alpha$-space 
and if $\,\alpha=\aleph_1\,$ then $\,X\,$ is a $P_\alpha$-space 
if and only if $\,X\,$ is a $P$-space. 
Clearly, if $\,X\,$ is a $P_\alpha$-space and $\,|X|<\alpha\,$ 
then $\,X\,$ is discrete.
(It is plain that if $\,X\,$ is a $P_\alpha$-space and $\,|X|=\alpha\,$ 
and $\,\alpha\,$ is a singular cardinal then $\,X\,$ is discrete.)
\smallskip
For an infinite cardinal $\,\alpha\,$ 
let us call a Hausdorff space {\it $\alpha$-normal}
when it is completely normal and every closed set 
is an intersection of at most $\,\alpha\,$ open sets.
So a Hausdorff space is perfectly normal 
if and only if it is $\aleph_0$-normal.
It is dull to consider 
{\it perfectly normal} $P$-spaces because, trivially, 
a perfectly normal $P$-space must be discrete.
More generally, if $\,\mu<\alpha\,$ then every $\mu$-normal 
$P_{\alpha}$-space is discrete.            
However, the enumeration problem concerning  
{\it completely normal} $P$-spaces
and $\alpha$-normal 
$P_{\alpha}$-spaces is not trivial and can be solved 
under certain cardinal restrictions.
\medskip
As usual, $\,\kappa^+\,$ denotes the smallest cardinal greater than 
$\,\kappa\,$, whence $\,\kappa^+\leq 2^\kappa\,$
and $\,\aleph_1=(\aleph_0)^+\,$.
Furthermore, for arbitrary $\,\kappa,\mu\,$ the cardinal number
$\,\kappa^{<\mu}\,$ is defined as usual (see [4]).
Note that if $\,\mu\leq\kappa^+\,$ then 
$\;\kappa^{<\mu}\,=\,|\{\,T\;|\;T\subset S\;\land\;|T|<\mu\,\}|\;$
whenever $\,S\,$ is a set of size $\,\kappa\,$.
In particular, $\,\kappa^{<\aleph_0}=\kappa\,$ 
and $\,\kappa^{<\aleph_1}=\kappa^{\aleph_0}\,$ 
for every $\,\kappa\,$.
Naturally, if $\,\mu=\kappa^+\,$ then $\,\kappa^{<\mu}=2^\kappa\,$.
Consequently, if $\,\mu>\kappa\,$ then $\,\kappa<\kappa^{<\mu}\,$.
(If $\,\kappa\,$ is a cardinal number of cofinality smaller than
$\,\mu^{++}\,$ then 
$\,\kappa<\kappa^{<\mu}\,$ due to K\"onig's Theorem [4, 5.14].)
On the other hand, for every $\,\mu\,$
the cardinals $\,\kappa\,$ 
satisfying $\,\kappa^{<\mu}=\kappa\,$ form a proper class $\,{\cal K}_\mu\,$
such that $\,2^\theta\in{\cal K}_\mu\,$ for every cardinal
$\,\theta\,$ with $\,\theta^+\geq\mu\,$
and if $\,\kappa\in{\cal K}_\mu\,$ 
then the cardinal successor $\,\kappa^+\,$
of $\,\kappa\,$ also lies in $\,{\cal K}_\mu\,$
due to the Hausdorff formula [4, (5.22)].
In particular, the cardinals 
$\;{\bf c},{\bf c}^+,{\bf c}^{++},...\;$ lie in $\,{\cal K}_\mu\,$
for $\,\mu=\aleph_1\,$. Furthermore, if $\,\kappa^{<\mu}=\kappa\leq \lambda\,$
and there are only finitely many cardinals $\,\theta\,$ 
with $\,\kappa\leq\theta\leq\lambda\,$ then $\,\lambda^{<\mu}=\lambda\,$.
(Note, again, that $\,\kappa^{<\alpha}=\kappa\,$
implies $\,\alpha\leq \kappa\,$.) 
\medskip
{\bf Theorem 8.} {\it Let $\,\alpha\,$ be an uncountable cardinal.
Assume $\,\kappa=\kappa^{<\alpha}\,$
and $\;\kappa\leq \lambda\leq 2^\kappa\;$ 
and $\,\lambda^{<\alpha}=\lambda\leq 2^\mu<2^\lambda\,$ 
for some $\,\mu\leq\kappa\,$ 
with $\,\mu^{<\alpha}=\mu\,$. 
Then there exist $\,2^{\lambda}\,$ mutually non-homeomorphic 
scattered, strongly zero-dimensional, 
hereditarily paracompact, $\alpha$-normal $P_\alpha$-spaces of size 
$\,\kappa\,$ and weight $\,\lambda\,$.
In particular, for every $\,\kappa\,$ with 
$\,\kappa=\kappa^{\aleph_0}\,$ there exist precisely 
$\,2^{2^\kappa}\,$ mutually non-homeomorphic paracompact P-spaces of size 
$\,\kappa\,$ and weight $\,2^\kappa\,$ up to homeomorphism.}
\medskip\smallskip                                
As usual (see [1] and [4]), a filter $\,{\cal F}\,$ is $\kappa$-{\it complete}
if and only if $\;\bigcap{\cal A}\in{\cal F}\;$ for every 
$\,{\cal A}\subset{\cal F}\,$ with $\,0\not=|{\cal A}|<\kappa\,$.
Trivially, every filter is $\aleph_0$-complete.
Obviously, an $\omega$-free filter is not $\kappa$-complete
for any $\,\kappa>\aleph_0\,$.
Let us call a filter $\,{\cal F}\,$ {\it $\kappa$-free} if and only if 
$\;\bigcap{\cal A}=\emptyset\;$ for some 
$\;{\cal A}\subset{\cal F}\;$ with $\;|{\cal A}|\leq\kappa\,$.
(So a filter is $\omega$-free if and only if it is $\aleph_0$-free.)
Clearly, if the topology of an almost discrete 
space $\,X\,$ is the single filter
topology defined with a free filter $\,{\cal F}\,$
then for every infinite cardinal $\,\alpha\,$
the (completely normal) space $\,X\,$ is $\alpha$-normal
if and only if $\,{\cal F}\,$ is $\alpha$-free,
and $\,X\,$ is a $P_\alpha$-space if and only if 
$\,{\cal F}\,$ is $\alpha$-complete.
Therefore, in view of the following counterpart
of Proposition 3, Theorem 8 can be easily proved 
by simply adopting the proof of the case $\,2^\lambda>2^\mu\,$ 
in Theorem~1 line by line while replacing the property
$\omega$-{\it free} with 
{\it $\alpha$-complete and $\alpha$-free} throughout.
\medskip
{\bf Proposition 7.} {\it If $\;|Y|=\kappa=\kappa^{<\mu}\;$ 
and $\;\kappa\leq\lambda=\lambda^{<\mu}\leq 2^\kappa\;$ 
then there exist $\,2^{\lambda}\,$ $\mu$-complete, 
$\mu$-free filters 
$\,{\cal F}\,$ on $\,Y\,$ such that $\;\chi({\cal F})=\lambda\,$.}
\medskip\smallskip
For the proof of Proposition 7 we need a lemma. 
This lemma also guarantees the existence of the family $\,{\cal A}_\omega\,$
in the proof of Proposition 3 since $\,\kappa^{<\mu}=\kappa\,$ for 
$\,\mu=\aleph_0\,$.
\medskip
{\bf Lemma 8.} {\it Let $\,Y\,$ be an infinite set 
of size $\,\kappa\,$ and assume  $\,\kappa^{<\mu}=\kappa\,$.
Then there exists a family
$\,{\cal A}\,$ of subsets 
of $\,Y\,$ such that $\,|{\cal A}|=2^\kappa\,$
and $\,{\cal A}\,$ has a subfamily $\,{\cal K}\,$ 
of size $\,\mu\,$ with $\,\bigcap{\cal K}=\emptyset\,$ 
and if $\;{\cal D},{\cal E}\not=\emptyset\;$ are 
disjoint subfamilies of $\,{\cal A}\,$ of size smaller than $\,\mu\,$
then $\,\bigcap{\cal D}\,$ is not a subset of $\,\bigcup {\cal E}\,$.}
\medskip
{\it Proof.} For an infinite set $\,S\,$ put 
$\;{\cal P}_\mu(S)\,:=\,\{\,T\;|\;T\subset S\;\land\; |T|<\mu\,\}\,$.
Let $\,Y\,$ be a set of size $\,\kappa\,$ and assume 
$\,\kappa^{<\mu}=\kappa\,$, whence $\,\kappa\geq\mu\,$.
Choose any set $\,X\,$ of size $\,\kappa\,$. 
Then $\;|{\cal P}_\mu(X)|=\kappa^{<\mu}=\kappa\;$ 
and hence 
$\;|{\cal P}_\mu({\cal P}_\mu(X))|=\kappa^{<\mu}=\kappa\,$. 
Therefore we may identify $\,Y\,$ with the set 
$\;{\cal P}_\mu(X)\times {\cal P}_\mu({\cal P}_\mu(X))\,$.
Now for $\;Y\;:=\;{\cal P}_\mu(X)\times 
{\cal P}_\mu({\cal P}_\mu(X))\;$ put 
\smallskip
\centerline{$\,A[S]\;:=\;\{\,(H,{\cal H})\in Y\;\;|\;\;
\emptyset\,\not=\,H\cap S\,\in\,{\cal H}\,\}\,$}
\smallskip
whenever $\;S\subset X\,$.
Clearly, $\,A[S]=\emptyset\,$ if and only if $\,S=\emptyset\,$.
We observe that $\;A[S_1]\not=A[S_2]\;$ whenever $\;S_1,S_2\subset X\;$
are distinct. Indeed, if $\,S_1,S_2\,$ are subsets of $\,X\,$ and 
$\;s\in S_1\setminus S_2\;$ then 
$\;(\{s\},\{\{s\}\})\,\in\,A[S_1]\setminus A[S_2]\,$.
Put $\;{\cal A}\;:=\;\{\,A[S]\;|\;\emptyset\not=S\subset X\,\}\,$.
Then $\,|{\cal A}|=2^\kappa\,$ and we 
claim that $\,{\cal A}\,$ is a family as desired.
\smallskip
For $\;0\not=|I\times J|<\mu\;$ let 
$\,\{\,S_i\;|\;i\in I\}\,$
and  $\,\{\,T_j\;|\;j\in J\}\,$ be disjoint families
of nonempty subsets of $\,X\,$. Choose 
$\;a_{i,j}\,\in\,(S_i\setminus T_j)\cup(T_j\setminus S_i)\;$
for every $\;(i,j)\,\in\,I\times J\;$
and $\,b_i\in S_i\,$ for every $\,i\in I\,$
and put $\;V\,:=\,\{\,a_{i,j}\;|\;i\in I,\,j\in J\,\}\cup
\{\,b_i\;|\;i\in I\,\}\,$.
Then $\,|V|<\mu\,$ and 
$\;\emptyset\not=V\cap S_i\not=V\cap T_j\;$ 
whenever $\,i\in I\,$ and $\,j\in J\,$. Hence
the pair $\;(V,\,\{\,V\cap S_i\;|\;i\in I\,\})\;$ 
lies in $\;\bigcap_{i\in I}A[S_i]\;$ but not in $\;\bigcup_{j\in J}A[T_j]\,$.
Finally, since $\,|H|<\mu\,$ whenever $\,(H,{\cal H})\in A[S]\,$, 
if $\,{\cal K}\,$ is any subfamily 
of $\;\{\,A[\{x\}]\;|\;x\in X\,\}\;$ with $\,|{\cal K}|=\mu\,$
then  $\,\bigcap{\cal K}=\emptyset\,$, {\it q.e.d.}
\medskip\smallskip
{\it Remark.} The previous proof is very similar to 
Hausdorff's classic construction of independent sets
as carried out in the proof of [4, 7.7].
However, by Hausdorff (and in [4, 7.7])
only the special case $\,\mu=\aleph_0\,$ is considered
and, unfortunately, from Hausdorff's construction one cannot   
obtain $\omega$-free resp.~$\alpha$-free
filters in a natural way. In order to accomplish this we have modified 
the proof of [4, 7.7] in a subtle but crucial way
by including the condition $\;\emptyset\not=H\cap S\;$
in our definition of $\,A[S]\,$. This condition 
guarantees that $\,{\cal A}\,$ has a subfamily 
$\,{\cal K}\,$ as desired and hence that the family $\,{\cal A}_\omega\,$
in the proof of Proposition 3 actually exists.
\bigskip                         
Now in order to prove Proposition 7
let $\,{\cal A}\,$ and $\,{\cal K}\,$ be families
as in Lemma 8. For every family $\,{\cal H}\,$ with 
$\,{\cal K}\subset{\cal H}\subset{\cal A}\,$
and $\,|{\cal H}|=\lambda\,\,$ put
\smallskip
\centerline{${\cal B}_{\cal H}\;:=\;
\{\,\bigcap{\cal G}\;
|\;\emptyset\not={\cal G}\subset{\cal H}
\;\land\;|{\cal G}|<\mu\,\}\,$.}
\smallskip
Then $\,\emptyset\not\in{\cal B}_{\cal H}\,$
and thus $\,{\cal B}_{\cal H}\,$ is a filter base for a $\mu$-complete
filter $\,{\cal F}[{\cal H}]\,$. 
Since $\,{\cal K}\subset{\cal F}[{\cal H}]\,$, the
filter $\,{\cal F}[{\cal H}]\,$ is $\mu$-free.
Since $\,\lambda^{<\mu}=\lambda\,$, we have 
$\;|{\cal B}_{\cal H}|=\chi({\cal F}[{\cal H}])=\lambda\;$
by exactly the same arguments as in the proof of Proposition 3.
\bigskip
{\it Remark.} Since for no cardinal $\,\kappa>\aleph_0\,$ 
the existence of a $\kappa$-complete {\it ultrafilter}
is provable in ZFC (see [4]), in Theorem 8 we cannot include the property 
{\it extremally disconnected}.
While Theorem~8 modifies Theorem~1 for $P$-spaces,
there is no pendant of Theorem 2 for $P$-spaces
because an infinite $P$-space is clearly not pathwise connected
and, moreover, every regular $P$-space $\,X\,$ is zero-dimensional.  
(If $\,x\in U\,$ where $\,U\subset X\,$ is open
then choose open neighborhoods $\,U_n\,$ of $\,x\,$ such that 
$\;U\supset \overline{U_n}\supset U_n\supset\overline{U_{n+1}}\supset U_{n+1}
\;$ for every $\,n\in{\Bbb N}\,$. 
Then 
$\;V\,:=\,\bigcap_{n=1}^\infty U_n=\bigcap_{n=1}^\infty \overline{U_n}\;$
is an open-closed neighborhood of $\,x\,$ and $\,V\subset U\,$.)
\bigskip\bigskip
{\bff References}
\medskip
[1] Comfort, W.W., and Negrepontis, S.: {\it The Theory of Ultrafilters.}
Springer 1974.
\smallskip
[2] Engelking, R.: {\it General Topology, revised and completed edition.}
Heldermann 1989. 
\smallskip
[3] Hodel, R.E.: {\it The number of metrizable spaces.}
Fund.~Math.~{\bf 115} (1983),
127-141.
\smallskip
[4] Jech, T.: {\it Set Theory.} 3rd ed. Springer 2002.
\smallskip
[5] Kuba, G.: {\it Counting topologies.}
Elemente d.~Math. {\bf 66} (2011), 56-62.
\smallskip
[6] Kuba, G.: {\it Counting metric spaces.}
Arch.~Math. {\bf 97} (2011), 569-578.
\smallskip
[7] Kuba, G.: {\it Counting linearly ordered spaces.}
Colloq.~Math. {\bf 135} (2014), 1-14.
\smallskip
[8] Kuba, G.: {\it On the variety of Euclidean point sets.}
Internat.~Math.~News~{\bf 228} (2015), 

\rightline{23-32.}

[9] Kuba, G.: {\it Counting ultrametric spaces.}
Colloq.~Math. {\bf 152} (2018), 217-234.
\smallskip
[10] Kulpa, W.: {\it On the existence of maps having graphs connected 
and dense.} 

\rightline{Fund.~Math.~{\bf 76} (1972), 207-211.}

[11] F.W.~Lozier and R.H.~Marty, {\it The number of continua.}
Proc.~Amer.~Math.~Soc.~{\bf 40} 

\rightline{(1973), 271-273.}

[12] S.~Mazurkiewicz and W.~Sierpi\'nski, {\it Contribution \`a la topologie
des ensembles d\'enom-} 

\rightline{{\it brables.} Fund.~Math.~{\bf 1} (1920), 17-27.}

[13] Steen, L.A., and Seebach Jr., J.A.:
{\it Counterexamples in Topology.} 
Dover 1995.
\bigskip\bigskip\bigskip
Gerald Kuba 
\smallskip
Institute of Mathematics, 

University of Natural Resources and Life Sciences, 1180 Wien, Austria. 
\smallskip
{\sl E-mail:} {\tt gerald.kuba@boku.ac.at}
\bigskip\medskip\bigskip\medskip\hrule\bigskip\medskip\bigskip
This paper is published under the less humorous title 
{\it Counting spaces of excessive weights}
in {\sl Matemati\v cki Vesnik}. 
(A revision of the title has been required
due to the very true fact that 
a mathematical definition of the property {\it overweight} 
cannot be found anywhere.)

\vfill\eject
\centerline{\bffg APPENDIX on Compact Hausdorff Spaces} 
\bigskip
In Section 9 we refer to the following important theorem from [7]. 
\medskip
{\bf Theorem A.} {\it For every $\,\kappa>\aleph_0\,$
there exist $\,2^\kappa\,$ mutually non-homeomorphic 
scattered, compact, linearly ordered spaces of size $\,\kappa\,$.}
\medskip
A proof of Theorem A has also been put on the {\bf arXiv}, see
\smallskip
\centerline{Kuba, G: {\it Scattered and paracomapct order topologies.} 
arXiv 2005.09368v1 (2020)} 
\smallskip
Notice that every linearly ordered space is {\it completely normal}
and that (cf.~[1, 2.24]) every totally disconnected 
(particularly, every scatterd) compact Hausdorff space is 
{\it strongly zero-dimensional.} Notice also   that
the size of a scattered, compact Hausdorff space
always coincides with its weight.
In particular, $\,2^\kappa\,$ in Theorem A is maximal 
due to (1.2). As already pointed out in Section 9, 
in the excluded case $\,\kappa=\aleph_0\,$ 
the statement in Theorem A would be true if and only if 
the Continuum Hypothesis were true. 
\medskip
In the following we 
derive two consequences of Theorem A worth mentioning. 
(Notice that always $\,\kappa^\kappa=2^\kappa\,$
and that if $\,\kappa\,$ is singular then 
$\,\kappa^\theta>\kappa\,$ for some infinite cardinal $\,\theta<\kappa\,$.)
\bigskip
{\bf Corollary A.} {\it Let $\,\kappa,\theta\,$ be cardinals
where either $\,\theta=1\,$ and $\,\kappa>\aleph_0\,$,
or $\,\aleph_0\leq \theta\leq \kappa\,$.
Then there exist precisely $\,2^\kappa\,$ 
totally disconnected, compact
Hausdorff spaces of weight $\,\kappa\,$ and 
size $\,\kappa^\theta\,$ up to homeomorphism.} 
\medskip
{\it Proof.} Despite (1.2) the cardinality $\,2^\kappa\,$ 
is {\it maximal} due to the estimate $\,|{\cal F}_2|\leq 2^\kappa\,$ 
in the remark in Section 10.
Thus it is enough to track down 
$\,2^\kappa\,$ mutually non-homeomorphic spaces as desired.  
The case that $\,\theta=1\,$ and $\,\kappa>\aleph_0\,$ is settled 
by Theorem A.
So we assume $\,\aleph_0\leq \theta\leq \kappa\,$ and, for the moment,
$\,\kappa>\aleph_0\,$ and we do not care whether
$\,\kappa^\theta=\kappa\,$ or $\,\kappa^\theta>\kappa\,$.
(Of course, in view of Theorem A there is nothing to show
if $\,\kappa^\theta=\kappa\,$ since $\,\kappa^\theta=\kappa\,$
implies $\,\kappa\geq c\,$.)
\smallskip
By virtue of Theorem A, 
we can define a family $\,{\cal H}_\kappa\,$ 
of scattered, compact Hausdorff spaces of size (and weight) $\,\kappa\,$
such that  $\,|{\cal H}_\kappa|=2^\kappa\,$
and distinct spaces in $\,{\cal H}_\kappa\,$ are never homeomorphic.
Fix one space $\,C\in {\cal H}_\kappa\,$
and consider the product space $\,C^\theta\,$.
Clearly, $\,C^\theta\,$ is a totally disconnected
compact Hausdorff space of size $\,\kappa^\theta\,$
and weight $\;\max\{w(C),\theta\}=\kappa\,.$
Furthermore, the space $\,C^\theta\,$ is obviously dense in itself.
For every $\,H\in{\cal H}_\kappa\,$
consider the topological sum $\,H+C^\theta\,$.
So $\,H+C^\theta\,$ is a totally disconnected compact Hausdorff space
of weight $\;\max\{|H|,w(C^\theta)\}=\kappa\;$
and size $\;\max\{|H|,|C^\theta|\}=\kappa^\theta\,$.
\smallskip
For distinct and hence non-homeomorphic spaces $\,H_1,H_2\in{\cal H}_\kappa\,$
the spaces $\,H_1+C^\theta\,$
and $\,H_2+C^\theta\,$ are not homeomorphic
because every $\,H\in {\cal H}_\kappa\,$ can be recovered
from $\,H+C^\theta\,$ in view of the obvious fact that $\,C^\theta\,$ 
is the perfect kernel of $\,H+C^\theta\,$. 
To finish the proof by settling 
the remaining case $\,\aleph_0=\theta=\kappa\,$
is left as a nice exercise. (Using Cantor derivatives it is not difficult 
to track down $\,{\bf c}\,$ mutually non-homeomorphic 
totally disconnected compact subspaces of the real line.)
\bigskip
By considering {\it cones} as in Section 6, it is straightforward 
to derive from Corollary A
the following noteworthy enumeration theorem about continua.
\medskip
{\bf Corollary B.} {\it Let $\,\kappa,\theta\,$ be cardinal numbers
where either $\,\aleph_0\leq \theta\leq \kappa\,$
or $\,\theta=1\,$. Then there exist precisely $\,2^\kappa\,$ 
pathwise connected, compact Hausdorff spaces of weight $\,\kappa\,$ and 
size $\,\max\{{\bf c},\kappa^\theta\}\,$.}     
\end